\documentclass[11pt]{amsart}

\usepackage[T1]{fontenc}
\usepackage[utf8]{inputenc}
\usepackage{lmodern}
\usepackage[a4paper,margin=27mm]{geometry}
\usepackage{amsmath,amssymb,amsthm,mathtools}
\usepackage{microtype}
\usepackage[colorlinks=true,linkcolor=blue,citecolor=blue,urlcolor=blue]{hyperref}

\newtheorem{theorem}{Theorem}[section]
\newtheorem{lemma}[theorem]{Lemma}
\newtheorem{proposition}[theorem]{Proposition}
\newtheorem{corollary}[theorem]{Corollary}
\theoremstyle{remark}

\numberwithin{equation}{section}

\DeclareMathOperator{\tr}{tr}
\DeclareMathOperator{\diag}{diag}
\DeclareMathOperator{\diam}{diam}
\DeclareMathOperator{\rank}{rank}
\DeclareMathOperator{\Sym}{Sym}
\DeclareMathOperator{\supp}{supp}
\DeclareMathOperator{\sgn}{sgn}

\newcommand{\R}{\mathbb R}
\newcommand{\G}{\overline{\Gamma}}
\newcommand{\norm}[1]{\left\lVert#1\right\rVert}

\newcommand{\ind}{\mathbf 1}
\newcommand{\B}{\mathrm B}
\newcommand{\p}{\partial}
\newcommand{\Mm}{\mathcal M^-_{1,\rho}}
\title[Sharp maximum estimates]
{Sharp maximum estimates for linear uniformly elliptic \\ equations with measurable coefficients}
\author{Shibing Chen, Neil S. Trudinger \&  Xu-Jia Wang}

\address{Shibing Chen}
\email{chenshibing1982@hotmail.com}

\address{Neil S. Trudinger, 
Mathematical Sciences Institute,
The Australian National University Canberra ACT 0200, Australia. }
\email{Neil.Trudinger@anu.edu.au}

\address{ Xu-Jia Wang, 
Institute for Theoretical Sciences, Westlake University, Hangzhou, 310030, China.}
\email{wangxujia@westlake.edu.cn}
\date{}
\thanks {This research was supported by Australian Research Council grant DP230100499 and  NSFC No.12225111.}
\begin{document}

\begin{abstract}
In this paper we obtain sharp exponents in the Aleksandrov-Bakelman-Pucci 
maximum estimate for linear uniformly elliptic equations with measurable coefficients,
originally conjectured by Pucci in 1966 \cite{Pucci}, together with corresponding sharp results under the more natural quantification of
uniform ellipticity in more than two dimensions.   The proofs combine a Hessian integral inequality obtained from boundary
Morse inequalities with an obstacle argument. The estimates also hold on arbitrary open sets, of finite volume,
 without regularity assumptions on their boundaries.
\end{abstract}

\maketitle

\baselineskip13.8pt
\parskip3pt

\section{Introduction}
Let $\Omega$ be a bounded domain in $\R^n$, $n\geq 2$.
Let $A=\{a_{ij}(x)\}$ be a measurable, symmetric, positive definite matrix.
Consider the elliptic operator in nondivergence form,
\begin{equation}\label{Operat}
 L_Au=A:D^2u=: \sum_{i,j=1}^n a_{ij}(x)D_{ij}u,\ \ \ x\in\Omega.
\end{equation}
We seek the precise range of $p$ such that the maximum estimate
\begin{equation}\label{eqmaximum}
 \sup_\Omega u\leq\sup_{\partial\Omega}u^+
 +C(\diam\Omega)^{2-n/p}\norm{(A:D^2u)^-}_{L^p(\Omega)}
\end{equation}
holds for functions $u$ in the class
\begin{equation}\label{eqfunction-class}
 u\in W^{2,p}_{\mathrm{loc}}(\Omega)\cap C(\overline\Omega) 
 \ \ \text{with}\ \  (A:D^2u)^-\in L^p(\Omega),
\end{equation} 
where $u^+=\max\{u, 0\}$ and $u^-=\max\{-u, 0\}$.

When $p=n$, \eqref{eqmaximum} is the famous 
Aleksandrov--Bakelman--Pucci maximum principle for linear uniformly elliptic equations; 
see \cite{Aleksandrov, Bakelman, Pucci66b, GT}. But the exponent $p=n$ is not optimal for particular operators.
For the Laplacian, the optimal range is $p>n/2$. 
For a general linear uniformly elliptic operator, 
the optimal exponent depends on an ellipticity constant. 
A question posed by Pucci is  how the optimal
exponent depends on the ellipticity condition \cite{Pucci,Pucci66b}.
For the normalized class
\begin{equation}\label{Cond-1}
 \tr A=1,\qquad A\geq\delta I,\qquad 0<\delta\leq1/n,
\end{equation}
Pucci \cite{Pucci} conjectured that the range is 
\begin{equation}\label{Range-1}
p>P_{\mathrm{tr}}(n,\delta):=\max\{n/2,n[1-(n-1)\delta]\} .
\end{equation}
The two lower bounds on the right hand side arise from different radial examples
which show that the range \eqref{Range-1} is sharp for $n>2$ \cite{Pucci}.

In dimension $n=2$, this conjecture was proved in \cite{AIM06}.
In higher dimensions, Fabes and Stroock \cite{FS84} obtained an
improvement below exponent $n$. Kuo and Trudinger \cite{KT07} proved
the maximum estimate under conditions expressed through the dual
G\aa rding cones, giving integer-exponent cases of the conjecture.
Trudinger \cite{T20} refined these results and the planar estimates.
For equations with first-order terms, Trudinger \cite{T26} obtained the
maximum estimate with improved integrability assumptions on the drift 
under dual-cone conditions. In particular in the uniformly elliptic case, under the normalisation
\eqref{Cond-1},  the drift coefficient $b$ need only 
lie in $L^q$ for some $q>n$, as in the Laplacian case. Aleksandrov had originally
 proved the case $p=n$ with $q=n$. Krylov \cite{Krylov 2021, Krylov 2021-2}  
 extended this to the uniformly elliptic cases where $p$ is sufficiently close to $n$, 
 and Jung and Woo \cite{JW25} recently gave an analytic
proof of such estimates.  See also Dong and Kitano \cite{DK25} for derivative
 estimates and critical exponents for fully nonlinear equations.

In dimension $n=2$, the Pucci quantification of uniform ellipticity through trace is equivalent to the
 standard one in terms of the maximum eigenvalue.  However, as already remarked in \cite{T20}, 
 Pucci conjecture alone is not the right one in higher dimensions. In this paper, as well as providing
  a proof of the full Pucci conjecture , we prove a further sharp result, using the maximum eigenvalue
 instead of the trace. 

First we state the Pucci conjecture result. 

\begin{theorem}\label{Thm1.1}
Let $n\geq3$ and assume \eqref{Cond-1}.  
For every $p>P_{\mathrm{tr}}(n,\delta)$, estimate \eqref{eqmaximum}
holds for all functions in the class \eqref{eqfunction-class}, 
where the constant depends only on $n,p$ and $\delta$. 
\end{theorem}

For the standard formulation of uniform ellipticity, we consider 
the  the second class of matrices satisfying
\begin{equation}\label{Cond-2}
 \lambda I\leq A(x)\leq\Lambda I,
\end{equation}
where $0<\lambda\leq\Lambda<\infty$.
Denote $\rho=\Lambda/\lambda\geq1$.
For $2\leq d\leq n$ and $1\leq m<d$, define
\[
 a_{d,m}=\frac{d-m}{m\rho},\qquad
 \lambda_j=\frac{d-j}{1+a_{d,m}},\qquad
 \mu_j=\frac{j a_{d,m}}{1+a_{d,m}}.
\]
Let $K_{d,m,\rho}$ be the $m\times m$ matrix given by
\begin{equation}\label{eqK-intro}
K_{d,m,\rho}=\begin{pmatrix}
\lambda_0+\mu_0, & -\lambda_0, &0, &\cdots, & 0\\
-\mu_1, & \lambda_1+\mu_1, &-\lambda_1, &\cdots, & 0\\
0, & -\mu_2, &\lambda_2+\mu_2, &\cdots, & 0\\
\cdot & \cdot & \cdot &\cdots & \cdot\\
0, & 0, &\cdots, &-\mu_{m-1}, & \lambda_{m-1}+\mu_{m-1}\\
\end{pmatrix} \end{equation}
Denote its smallest eigenvalue by
$\kappa_{d,m}(\rho)$. One has
\[
 \kappa_{n,1}(\rho)=\frac{n\rho}{\rho+n-1}.
\]
Set
\begin{equation}\label{eqbox-threshold}
 P_n(\rho)=\max\Big\{\frac n2, \ \ 
       \max_{\substack{2\leq d\leq n, 1\leq m<d}}
         \  \kappa_{d,m}(\rho)\Big\}.
\end{equation}
If $n=2$, one can verify directly that $P_2(\rho)=P_{\mathrm{tr}}(2,1/(1+\rho))$.
We will show that
\begin{equation} \label{Range-2}
 n/2\leq P_n(\rho)<n. 
\end{equation}

For the standard formulation of uniform ellipticity \eqref{Cond-2}, we have the following sharp maximum estimate.
\begin{theorem}\label{Thm1.2}
Let $n\geq 3$ and assume \eqref{Cond-2}.  
For every $p>P_n(\rho)$, estimate \eqref{eqmaximum}
holds for all functions in the class \eqref{eqfunction-class},  
where the constant $C$ depends only on $n,p,\lambda,\Lambda$. 
\end{theorem}

The exponents $P_{\mathrm{tr}}(n,\delta)$ in Theorem \ref{Thm1.1} and $P_n(\rho)$ in Theorem \ref{Thm1.2} are optimal. 
When $1\le p\leq P_{\mathrm{tr}}(n,\delta)$ in Theorem \ref{Thm1.1} or $1<p< P_n(\rho)$ in Theorem \ref{Thm1.2},
estimate \eqref{eqmaximum} does not hold, even for smooth functions in the unit ball, vanishing on the boundary.

A trace-one matrix with $A\geq\delta I$ also satisfies
$A\leq[1-(n-1)\delta]I$. Conversely, if $I\leq A\leq\rho I$, then
$ \frac A{\tr A}\geq\frac I{1+(n-1)\rho}.$
The ranges of the exponent $p$ in the above two theorems partially overlap, 
but neither is completely contained in the other.

We also remark here that from the proof of Theorem 1.1 in \cite{T20}, we can extend Theorems 1.1 and 1.2 to arbitrary domains $\Omega$, with the dependence on ($\diam\Omega$) replaced by that on $|\Omega| ^{1/n}$, thereby providing an extension to unbounded domains.  Corresponding extensions of the local maximum principle, weak and strong Harnack inequalities and Holder estimates of Trudinger, and Krylov -Safonov, as presented for example in  \cite{GT}, also follow.  Concerning extensions to operators with first order terms we also have from the proof of the case $k \le n/2$ in Theorem 1.1 of \cite{T26}, that Theorems 1.1 and 1.2 extend to such operators,

\begin{equation}\label{Operat-b}
 Lu=A:D^2u+b\cdot Du.
\end{equation}
with drift coefficient $b$ in $L^q(\Omega)$, with $q=2p$, which still does not extend the integer case in \cite{T26} for $p > n/2$. More generally we need only assume the quantity $A^{-1}b\cdot b$ lies in  $L^p(\Omega)$.

The following is an outline of the proof of Theorems \ref{Thm1.1} and \ref{Thm1.2}. 

\begin{enumerate}
\item Choose a noninteger exponent $q$ strictly between the proposed
threshold and $\min\{p,n\}$, and put $k=\lfloor q\rfloor$, the largest integer under $q$.
The integer $k$ selects the Hessian cone $\Gamma_k$, while $q$ is the
integrability exponent that the estimate must retain. For
$H\in\G_k$, let $b$ be its $k$th largest eigenvalue. Subtract the
first $k+1$ terms of a characteristic polynomial and integrate the
remainder against $t^{q-n-1}\,dt$. This defines a continuous function
$Q_{n,k,q}(H)$, homogeneous of degree $q$.

\item Establish a nonnegative integral for $Q$, rather than a pointwise
sign. Apply the strong boundary Morse inequalities to
$t(|x|^2-R^2)/2-w-y\cdot x$ and integrate over $y$.
The signed area formula turns the critical-point counts into a truncated
determinant integral. Newton's boundary identities then express the
subtracted terms through the normal derivative of $w$ on the sphere.
Integration in $t$ gives
$\int_{B_R}Q_{n,k,q}(D^2w)\geq0$ for every admissible $w$ vanishing
on the boundary. The proof also covers a vanishing selected eigenvalue.

\item Prove two matrix bounds. The first,
$Q_{n,k,q}(H)\leq C\sigma_k(H)^{q/k}$, holds on the whole cone.
The second has the form
\[
 Q_{n,k,q}(H)\leq C(F(H)^+)^q-c|H|^q,
\]
where $F(H)$ is the infimum of $A:H$ over the coefficient class.
To obtain it, it is enough to prove $Q(H)<0$ when $H\ne0$ and
$F(H)\leq0$; compactness on $|H|=1$ supplies uniform constants.
For fixed trace, separate the smallest eigenvalue. The expression is
affine in that eigenvalue, and at the critical value it becomes a
weighted sum over deleted pairs. A symmetric compression represents the
whole sum. Interlacing and Newton's inequalities put its roots in the
cone needed for a fractional integral sign. Individual summands need
not have that sign.

\item For prescribed eigenvalue bounds, first average the negative
eigenvalues. Then maximise over positive eigenvalues with their sum
fixed. An extremum argument, including repeated selected eigenvalues,
reduces the problem to spectra of the form
$((-a)^m,1^{d-m})$. Zero eigenvalues account for dimensions $d<n$.
The elementary symmetric functions of these spectra are values of a
Krawtchouk polynomial. Its first zero is the smallest eigenvalue of
$K_{d,m,\rho}$. The condition $q>P_n(\rho)$ puts $q$ beyond every
relevant first zero and gives the strict matrix sign.

\item Use an upper obstacle to connect the matrix inequalities to the
original function. After subtracting its boundary maximum, assume that
$u\leq0$ on $\partial B_R$. Solve below $-u+\varepsilon$, with zero
boundary values and $\sigma_k(D^2z)^{1/k}=\tau$ away from contact.
At contact, $D^2z\leq-D^2u$, so
$F(D^2z)\leq(A:D^2u)^-$.
Use the coercive bound on the contact set and the whole-cone bound
elsewhere. Together with the nonnegative integral, they give
\[
 \norm{\sigma_k(D^2z)^{1/k}}_{L^q}
 \leq C\norm{(A:D^2u)^-}_{L^q}+C\tau R^{n/q}.
\]
The Hessian potential estimate bounds the height of $z$ by this norm.
Its radial integral converges precisely because $q>n/2$.
Since $u-\varepsilon\leq-z$, sending $\tau$ and $\varepsilon$ to zero
proves the ball estimate.

\item Raise $q$ to $p$ by H\"older's inequality. Smooth approximation
gives the Sobolev estimate on a ball. While we can use the Pucci extremal equation 
method in \cite{T20} to extend to arbitrary open bounded sets, we proceed somewhat
 differently here. On an arbitrary bounded open set, subtract the boundary maximum
  and a positive constant, take a smooth convex approximation to the positive part, 
  and extend by zero into a containing ball. The resulting function has compact 
  support inside the original set. Its Hessian has an extra positive semidefinite term,
which can only improve the one-sided forcing bound.

\item Prove necessity. Radial solutions with logarithmic cutoffs give
failure at $n/2$ and at the fixed-trace threshold. For the spectral
threshold, split diagonal matrices along symmetric rank-one directions,
stopping each branch when it has enough negative entries. The same
matrix $K_{d,m,\rho}$ governs the decay of the remaining mass.
Below its smallest eigenvalue, the forcing energy tends to zero.
At equality, a small change in the negative eigenvalues creates a
spectral gap; the accumulated forcing is then of order
$\eta^{p-1}$ and still tends to zero. Realise these finite laminates
by smooth Hessians and choose smooth coefficients. A convex truncation
and rescaling extend the examples from dimension $d$ to dimension $n$
without losing the critical exponent.

\end{enumerate}

To implement the above argument, we will draw on analytic techniques from works on Hessian equations,
including the admissible Hessian theory of \cite{CNS,TrudingerWang}, 
the obstacle regularity theorem of Jiao and Wang \cite{JW14}, 
and the potential estimate of Phuc and Verbitsky \cite{PV08};
as well as boundary Morse inequalities in the convention of Laudenbach \cite{Lau11}.
For the lower estimates, we will use the finite symmetric laminate realization of Astala, Faraco, and Székelyhidi \cite{AFS}, following the convex integration methods of Müller and Šverák \cite{MS} and the laminate constructions of Faraco \cite{Faraco}.

The arrangement of the paper is as follows. Section 2 collects the notation and basic properties of the Hessian cones. 
Sections 3 and 4 establish the Hessian integral inequality and the matrix estimates used in the obstacle argument of Section 5. 
Section 6 proves the matrix inequality for coefficients with fixed trace. 
Sections 7 and 8 treat general ellipticity bounds and identify the spectral threshold \(P_n(\rho)\). 
These results are combined in Section 9 to prove the upper estimates in Theorems 1.1 and 1.2. 
Section 10 gives the radial counterexamples, proving the sharpness of Theorem 1.1. 
Sections 11 and 12 construct the spectral counterexamples and extend them to higher dimensions, 
completing the sharpness proof for Theorem 1.2.

\underline {Acknowledgements}

The proofs of both Theorems \ref{Thm1.1} and \ref{Thm1.2} were discovered by ChatGPT more than three weeks ago. The authors have checked and rewritten such input and take responsibility for the full proof.

A different proof of Theorem \ref{Thm1.1} was found, independently by Hung and Wu, (In September, just after ours), and is presented in \cite{Hung-Wu}. While reviewing and organizing this work, we also learnt that the authors of \cite{LMMNSZ} also utilized AI to find a different proof.  However, the proof in this paper differs from those, and our Theorem \ref{Thm1.2} is not included in their findings.

 More recently the planar case has been treated differently from the approach in \cite{AIM06}  by  Memana, Trudinger and Yung \cite{MTY}  leading to the extension of the Pucci conjecture to the non-strict inequality $p\ge P_{\mathrm{tr}}(n,\delta)$ for dimension $n=2$ and exponent $p > 1$. The authors of this paper are grateful to the other authors of these forthcoming papers, for sharing their methods and results

\vskip20pt

\section{Preliminaries}
For $H\in\Sym(n)$, write its eigenvalues in non-decreasing order as
$\lambda_1(H)\leq\cdots\leq\lambda_n(H)$. Write $\sigma_j(H)$ for
the $j$th elementary symmetric function of $\lambda=(\lambda_1, \cdots, \lambda_n)$, 
with $\sigma_0=1$.
We use the cones
\[
 \Gamma_k=\{H:\sigma_j(H)>0\ (1\leq j\leq k)\},\qquad
 \G_k=\overline{\Gamma_k},
\]
and use the same notation for eigenvalue vectors $\lambda$. 
For convenience we denote $\G_0=\R^n$ for eigenvalue vectors and $\G_0=\Sym(n)$ for matrices, and
let  $|H|=(\tr H^2)^{1/2}$ be the Frobenius norm.
A $C^{1,1}$ function is $k$-admissible if its Hessian belongs to
$\G_k$ almost everywhere.

\begin{lemma}\label{lemcone}
Let $1\leq k\leq n$ and $\lambda\in\overline\Gamma_k$.
Assume that $\lambda_1\leq\cdots\leq\lambda_n$. Then 
\begin{itemize}
\item [(i)] $\sigma_k=0$ on $\partial\Gamma_k$ and  $\sigma_k^{1/k}(\lambda)$ is concave on $\Gamma_k$.
\item [(ii)]  $\frac{\p}{\p \lambda_j}\sigma_k(\lambda)\ge 0\ \ \forall\ j=1, \cdots, n$.
\item [(iii)] Let $m=n-k+1$. Then $b :=\lambda_m \ge 0$, and
\begin{equation}\label{eqsmall-eigenvalues}
{\begin{split}
   \sum_{i=1}^{m}\lambda_i & \geq0,\\ 
 -(m-1)b  \leq\lambda_i & \leq b\quad \ (i\leq m).
 \end{split}}
\end{equation}
\item [(iv)] We have $b>0$ for $\lambda\in\Gamma_k$.
\item [(v)]
If $b=0$, the matrix is positive semidefinite and has rank at most $k-1$. 
\end{itemize}
\end{lemma}

\begin{proof} Properties (i) and (ii) are well known. For (iii), 
delete the largest $k-1$ eigenvalues. The remaining list belongs to
$\G_1$, so its sum is nonnegative. Its largest entry is $b$, which
must therefore be nonnegative. For each remaining entry, subtracting
the sum of the other $m-1$ entries gives the lower bound in
\eqref{eqsmall-eigenvalues}. In the open cone the deleted list is
in $\Gamma_1$, and its sum is positive, so $b>0$.
If $b=0$, the lowest $m$ entries are nonpositive with nonnegative
sum, hence all vanish. The remaining $k-1$ entries are nonnegative.
\end{proof}

If a negative entry $-a$ is deleted from a vector in $\G_k$, the
remaining list $\xi$ belongs to $\Gamma_k$, provided it has length at
least $k$. Indeed,
\[
 \sigma_j(\xi)=\sigma_j(H)+a\sigma_{j-1}(\xi)>0
 \quad(1\leq j\leq k),
\]
by induction from $\sigma_0=1$.
The two lower extremal operators are
\begin{equation}\label{eqextremals}
 \begin{aligned}
 F_\delta(H)&=\delta\tr H+(1-n\delta)\lambda_1(H),\\
 \Mm(H)&=\sum_{\lambda_i(H)>0}\lambda_i(H)
          +\rho\sum_{\lambda_i(H)<0}\lambda_i(H).
\end{aligned}
\end{equation}
They are the infima of $A:H$ over $\{A\geq\delta I:\tr A=1\}$
and $\{I\leq A\leq\rho I\}$, respectively. Both are continuous and
positively homogeneous of degree one. Finally, write
\[
 \B(s,t)=\int_0^1r^{s-1}(1-r)^{t-1}\,dr\qquad(s,t>0).
\]

\vskip20pt

\section{Boundary Morse inequalities and Hessian integrals}
 We next obtain an integral inequality without imposing any cone condition.
For $0\leq r\leq n-2$, define
\[
 G_r(M)=\bigg\{{\begin{aligned}
 &\det M\ \ \text{if}\ \ \lambda_{r+1}(M)\geq0,\\
 &0\hskip31pt \text{otherwise.}
 \end{aligned}}
\]
Denote by $\nu$ the unit outer normal.

\begin{lemma}\label{lemmorse}
If $U\in C^{1,1}(\overline{B_R})$ and $U=0$ on $\partial B_R$, then
\begin{equation}\label{eqmorse}
 (-1)^r\left[  \int_{B_R}G_r(D^2U)\,dx
 -\frac{1}{nR^{n-1}}\int_{\partial B_R}(U_\nu^+)^n\,dS \right]
 \geq0,\qquad 0\leq r\leq n-2.
\end{equation}
\end{lemma}

If $U$ is uniformly convex, 
then $G_r(D^2 U)$ is independent of $r$, and one has
$$ \int_{B_R}G_r(D^2U)\,dx
 = \frac{1}{nR^{n-1}}\int_{\partial B_R}(U_\nu^+)^n\,dS .$$

\begin{proof}
Assume first that $U$ is smooth on $\overline{B_R}$.
For $y\in\R^n$, set
\[
 f_y(x)=U(x)-y\cdot x.
\]
The critical points of $f_y$ satisfy $DU(x)=y$, and their
Hessians are $D^2U(x)$. For almost every $y\ne0$, $y$ is a
regular value of $DU$ and $y\notin DU(\partial B_R)$.
The corresponding critical points are therefore nondegenerate
and form a finite set in the interior.
Let $m_j(y)$ denote the number of these points whose Hessian
has exactly $j$ negative eigenvalues.

The boundary critical points are particularly simple.
Since $U=0$ on the sphere,
\[
 f_y|_{\partial B_R}=-y\cdot x.
\]
This function has two critical points:
\[
 \begin{aligned}
 x_+(y)&=R\frac y{|y|},
 &&\text{a minimum of index }0,\\
 x_-(y)&=-R\frac y{|y|},
 &&\text{a maximum of index }n-1.
 \end{aligned}
\]
Here the boundary index is computed from the Hessian in the
$n-1$ tangential directions. At the minimum,
\[
 \partial_\nu f_y(x_+(y))
 =U_\nu(x_+(y))-|y|.
\]

We now apply the boundary Morse inequalities.
Let $c_j$ count the interior critical points of index $j$,
together with the boundary critical points of tangential
index $j$ for which $\partial_\nu f_y<0$.
The latter condition means that $f_y$ increases into the ball.
Near such a point, suitable local coordinates $(z,s)$,
with $s\geq0$ pointing inward, give
\[
 f_y=f_y(x_0)+s+q(z),
\]
where $q$ is a nondegenerate quadratic form of index $j$.
These boundary points contribute to changes in the topology
of the sublevel sets.

Since the ball is contractible, the polynomial form of the
boundary Morse inequalities
\cite[Theorem A and Corollary A]{Lau11} is
\[
 \sum_{j=0}^n c_jt^j-1
 =(1+t)\sum_{j\geq0}a_jt^j,
 \qquad a_j\geq0.
\]
Taking partial alternating sums gives
$\sum_{j=0}^r(-1)^jc_j-1=(-1)^ra_r$, and hence
\begin{equation}\label{eqmorse-ball-count}
 (-1)^r\Big(\sum_{j=0}^r(-1)^jc_j-1\Big)\geq0.
\end{equation}

For $r\leq n-2$, the boundary maximum has index greater than
$r$ and does not enter this sum. Thus, writing
\[
 N_0(y)=\ind_{\{U_\nu(Ry/|y|)<|y|\}},
\]
we have $c_0=m_0(y)+N_0(y)$ and $c_j=m_j(y)$ for
$1\leq j\leq r$. It follows that
\begin{equation}\label{eqmorse-y-count}
 (-1)^r\Big[
 \sum_{j=0}^r(-1)^jm_j(y)-(1-N_0(y))
 \Big]\geq0
 \quad\hbox{for almost every }y.
\end{equation}

To convert this count into an integral, let $E_j\subset B_R$
be the set where $D^2U$ is nonsingular and has exactly $j$
negative eigenvalues. On this set,
\[
 \det D^2U=(-1)^j|\det D^2U|.
\]
Moreover, $\lambda_{r+1}(D^2U)\geq0$ precisely when $j\leq r$.
Since $G_r$ vanishes at singular Hessians, the area formula gives
\begin{equation}\label{eqmorse-area-count}
 \begin{aligned}
 \int_{B_R}G_r(D^2U)\,dx
 &=\sum_{j=0}^r(-1)^j
       \int_{E_j}|\det D^2U|\,dx\\
 &=\int_{\R^n}\sum_{j=0}^r(-1)^jm_j(y)\,dy.
 \end{aligned}
\end{equation}
For the other term, write $y=s\theta$, with $s>0$ and
$\theta\in\mathbb S^{n-1}$. By the definition of $N_0$,
\[
 1-N_0(s\theta)=\ind_{\{s\leq U_\nu(R\theta)\}}.
\]
Integration in polar coordinates therefore yields
\begin{align*}
 \int_{\R^n}(1-N_0(y))\,dy
 &=\int_{\mathbb S^{n-1}}
     \int_0^{U_\nu(R\theta)^+}s^{n-1}\,ds\,d\theta\\
 &=\frac1{nR^{n-1}}
     \int_{\partial B_R}(U_\nu^+)^n\,dS.
\end{align*}
Integrating \eqref{eqmorse-y-count} and using these identities
proves \eqref{eqmorse} for smooth $U$.

Finally, let $U\in C^{1,1}(\overline{B_R})$.
Fix $s>n$, choose $g_i\in C^\infty(\overline{B_R})$
converging to $\Delta U$ in $L^s(B_R)$, and solve
\[
 \Delta U_i=g_i\quad\hbox{in }B_R,\qquad
 U_i=0\quad\hbox{on }\partial B_R.
\]
The solutions are smooth on the closed ball. The Dirichlet
estimate gives
\[
 \norm{U_i-U}_{W^{2,s}(B_R)}
 \leq C\norm{g_i-\Delta U}_{L^s(B_R)}
 \longrightarrow0.
\]
Consequently, $U_i\to U$ in $C^1(\overline{B_R})$ and
$D^2U_i\to D^2U$ in measure.

The function $G_r$ is continuous: its indicator changes only
where $\lambda_{r+1}=0$, and the determinant vanishes there.
Thus $G_r(D^2U_i)\to G_r(D^2U)$ in measure.
Also, $|G_r(H)|\leq C_n|H|^n$, so for every measurable
$E\subset B_R$,
\[
 \int_E|G_r(D^2U_i)|\,dx
 \leq C_n\norm{D^2U_i}_{L^s(B_R)}^n|E|^{1-n/s}
 \leq C|E|^{1-n/s}.
\]
Since $s>n$, this gives uniform integrability. Together with
convergence in measure, it implies convergence of the bulk
integrals. The boundary integrals converge by the uniform
convergence of $DU_i$. Passing to the limit proves
\eqref{eqmorse}.
\end{proof}

Fix $1\leq k\leq n-1$ and $k<p<k+1$. For $H\in\G_k$ put
$b=\lambda_{n-k+1}(H)$ and define
\begin{equation}\label{eqQ}
 Q_{n,k,p}(H)=
 \begin{cases}
 \displaystyle(-1)^k\sum_{j=0}^n
 \frac{(-1)^j\sigma_j(H)b^{p-j}}{p-j},&b>0,\\[2mm]
 0,&b=0.
 \end{cases}
\end{equation}
The selected eigenvalue $b$ is the $k$th largest eigenvalue.

\begin{proposition}\label{propintegral}
The function $Q_{n,k,p}$ is continuous on $\G_k$, is positively homogeneous
of degree $p$, and satisfies $|Q_{n,k,p}(H)|\leq C|H|^p$.
For every zero-boundary $k$-admissible
$w\in C^{1,1}(\overline{B_R})$,
\begin{equation}\label{eqQ-integral}
 \int_{B_R}Q_{n,k,p}(D^2w)\,dx
 \geq\frac{(n-1)!}{\left|\prod_{j=0}^n(p-j)\right|}
 R^{1-p}\int_{\partial B_R}w_\nu^p\,dS\geq0.
\end{equation}
\end{proposition}

\begin{proof}
By Lemma~\ref{lemcone}, the lowest $n-k+1$ eigenvalues have absolute
value at most $C_nb$. A product of $j\geq k$ eigenvalues contains at
least $j-k+1$ of these entries. Consequently
\[ {\begin{aligned}
 |\sigma_j(H)| & \leq C_n b^{j-k+1}|H|^{k-1}\quad(j\geq k), \\
  |\sigma_j(H)| &\leq C_n|H|^j\quad(j<k).
 \end{aligned}}
\]
After multiplication by $b^{p-j}$, these inequalities give the asserted
growth and continuity at $b=0$, because $p-k+1>0$.
Homogeneity and continuity where $b>0$ follow directly from the definition.

To prove the integral inequality, subtract the first $k+1$ terms of the
characteristic polynomial and set
\[
 L_t(H)=G_{k-1}(tI-H)-\sum_{j=0}^k(-1)^j\sigma_j(H)t^{n-j}.
\]
The indicator in $G_{k-1}(tI-H)$ switches at $t=b$. Therefore
\[
 L_t(H)=
 \begin{cases}
 -\displaystyle\sum_{j=0}^k(-1)^j\sigma_j(H)t^{n-j},&0<t<b,\\
 \displaystyle\sum_{j=k+1}^n(-1)^j\sigma_j(H)t^{n-j},&t>b.
 \end{cases}
\]
Integrating each monomial gives
\begin{equation}\label{eqmellin}
 (-1)^{k-1}\int_0^\infty t^{p-n-1}L_t(H)\,dt=Q_{n,k,p}(H).
\end{equation}
The first integral converges at zero since $p>k$, and the second at
infinity since $p<k+1$. The same eigenvalue bounds as above show that
its absolute integral is at most $C|H|^p$. If $b=0$, the rank assertion
in Lemma~\ref{lemcone} gives $L_t(H)=0$.

Let $q=w_\nu$ on the sphere. Since $\Delta w\geq0$ and $w=0$ on the
boundary, $w\leq0$ and $q\geq0$. We use the Newton tensor boundary identities \cite{Reilly}
\begin{equation}\label{eqnewton-boundary}
 \int_{B_R}\sigma_j(D^2w)\,dx
 =\frac1n\binom nj R^{1-j}\int_{\partial B_R}q^j\,dS,
 \qquad 1\leq j\leq n.
\end{equation}
For smooth $w$, the Newton tensor
\[
T_{j-1}^{ab}=\partial\sigma_j(D^2w)/\partial w_{ab}
\] 
is divergence free, and
\[
j\sigma_j(D^2w)=\partial_a(T_{j-1}^{ab}w_b) .
\]
On the sphere, $Dw=q\nu$ and the tangential Hessian is $(q/R)I$.
Hence 
\[
T_{j-1}^{\nu\nu}=\binom{n-1}{j-1}(q/R)^{j-1},
\] 
which gives \eqref{eqnewton-boundary}. 
The Poisson approximation used in
Lemma~\ref{lemmorse} extends this identity to $C^{1,1}$ functions.
For $j=0$, the corresponding formula is simply
$|B_R|=R|\partial B_R|/n$.

Apply Lemma~\ref{lemmorse} with $r=k-1$ to
$U=t(|x|^2-R^2)/2-w$. Subtracting the integrals in
\eqref{eqnewton-boundary} gives
\[
 (-1)^{k-1}\left[\int_{B_R}L_t(D^2w)\,dx -\int_{\partial B_R}\Psi_t(q)\,dS\right] \geq 0,
\]
where
\[
 \Psi_t(q)=\frac Rn\Big[(t-q/R)_+^n -\sum_{j=0}^k(-1)^j\binom nj(q/R)^jt^{n-j}\Big].
\]
Multiply by $t^{p-n-1}$ and integrate over $t>0$.
For $q>0$, split the scalar integral at $t=q/R$ and integrate its
monomials, exactly as in \eqref{eqmellin}. The identity
\[
 \sum_{j=0}^n\frac{(-1)^j\binom nj}{p-j}
 =\frac{(-1)^nn!}{\prod_{j=0}^n(p-j)}
\]
then gives
\[
 (-1)^{k-1}\int_0^\infty t^{p-n-1}\Psi_t(q)\,dt
 =\frac{(n-1)!R^{1-p}q^p}{\left|\prod_{j=0}^n(p-j)\right|}.
\]
For $q=0$, both sides are zero. Absolute convergence, with the preceding
growth bounds and boundedness of $D^2w$ and $q$, justifies Fubini's
theorem and proves \eqref{eqQ-integral}.
\end{proof}

\vskip20pt

\section{Matrix estimates}
The obstacle argument needs an upper bound on the whole Hessian cone.
We first establish two polynomial signs, then use them to prove the
bound and a fractional sign inequality used in both coefficient classes.
\begin{lemma}\label{lemtaylor}
Let $h_1\geq\cdots\geq h_d$ belong to $\G_l$, where $1\leq l\leq d$,
and suppose $c=h_l>0$. Set
\[
 \begin{aligned}
 D_h(t)&=\prod_{i=1}^d(1-th_i),\\
 E_l(h,t)&=D_h(t)-\sum_{j=0}^{l-1}(-1)^j\sigma_j(h)t^j,\\
 T_l(h,t) &=(-1)^{l-1}\sum_{j=0}^{l-1}(-1)^j\sigma_j(h)t^j.
 \end{aligned}
\]
Then
\[  \begin{aligned}
 (-1)^lE_l(h,t) &\geq 0\ \ \ \ \ (0\leq t\leq c^{-1}),\\
 T_l(h,t) &\geq (ct-1)^{l-1}\quad(t\geq c^{-1}).
  \end{aligned}
\]
\end{lemma}

\begin{proof}
Put $L=h_1\geq c$. For $0\leq s<L^{-1}$, every factor in $D_h(s)$ is
positive and
\[
 \frac{h_i}{1-sh_i}-h_i=\frac{s h_i^2}{1-sh_i}\geq0.
\]
The vector with entries $h_i/(1-sh_i)$ is therefore in $\G_l$.
Differentiating the product gives
\[
 (-1)^l D_h^{(l)}(s)
 =l!D_h(s)\sigma_l\left(\frac{h_1}{1-sh_1},\ldots,
                         \frac{h_d}{1-sh_d}\right)\geq0.
\]
Taylor's integral formula proves the first sign up to $t=L^{-1}$.
This proves its full range when $l=1$.

For $l\geq2$, delete $L$ and call the remaining vector $h'$. It belongs
to $\G_{l-1}$ and its $(l-1)$st largest entry is $c$. Expanding
$D_h(t)=(1-Lt)D_{h'}(t)$ gives
\begin{align*}
 (-1)^lE_l(h,t)
 &=(Lt-1)(-1)^{l-1}E_{l-1}(h',t)
   +\sigma_{l-1}(h')t^{l-1},\\
 T_l(h,t)&=(Lt-1)T_{l-1}(h',t)+\sigma_{l-1}(h')t^{l-1}.
\end{align*}
The first identity and induction extend the first sign over
$L^{-1}\leq t\leq c^{-1}$.
For the second sign, start with $T_1=1$. If $t\geq c^{-1}$, then
$Lt-1\geq ct-1\geq0$. The second identity, deletion, and induction
give $T_l(h,t)\geq(ct-1)^{l-1}$.
\end{proof}

\begin{proposition}\label{propupper}
For $1\leq k\leq n-1$ and $k<p<k+1$, there is a constant $C=C(n,k,p)$
such that
\[
 Q_{n,k,p}(H)\leq C\sigma_k(H)^{p/k}\qquad(H\in\G_k).
\]
\end{proposition}

\begin{proof}
We first show a strict sign when $\sigma_k(H)=0$ and
$b=\lambda_{n-k+1}(H)>0$. Use $E_k$ and $T_k$ from
Lemma~\ref{lemtaylor}, applied to the eigenvalues of $H$.
Since $\sigma_k(H)=0$, the remainder $E_k(H,t)$ is $O(t^{k+1})$
at zero. Termwise integration in \eqref{eqQ} gives
\begin{align*}
 Q_{n,k,p}(H)
 &=(-1)^k\sum_{j=0}^{k-1}
       \frac{(-1)^j\sigma_j(H)b^{p-j}}{p-j}
   -\int_0^{1/b}t^{-p-1}(-1)^kE_k(H,t)\,dt\\
 &=-\int_{1/b}^{\infty}t^{-p-1}T_k(H,t)\,dt
   -\int_0^{1/b}t^{-p-1}(-1)^kE_k(H,t)\,dt\\
 &\leq-\int_{1/b}^{\infty}t^{-p-1}(bt-1)^{k-1}\,dt
 =-b^p\B(p-k+1,k)<0.
\end{align*}
Here $\B(a,b)=\int_0^1s^{a-1}(1-s)^{b-1}\,ds$.
Both integrals converge because $k<p<k+1$.

Normalize $b=1$ and compactify the largest $k-1$ eigenvalues
by taking their reciprocals. Write
\[
 H=\diag(a_1,\ldots,a_m,L_1,\ldots,L_{k-1}),
 \]
 where
$m=n-k+1,\quad a_m=1,\quad a_i\leq1,\quad L_i\geq1.$
Lemma~\ref{lemcone} gives $a_i\geq-(m-1)$.
Set
\[
 P=\prod_{i=1}^{k-1}L_i\geq1,\qquad r_i=L_i^{-1},\qquad
 S=\frac{\sigma_k(H)}P,\qquad F=\frac{Q_{n,k,p}(H)}P.
\]
For $k=1$, the product is understood to be $1$.
Both $S$ and $F$ are polynomials in $(a,r)$: division by $P$ replaces
each omitted large eigenvalue by its reciprocal. Thus they extend
continuously to the compact closure $K$ of the admissible normalized
vectors in the box $a_i\in[-(m-1),1]$, $r_i\in[0,1]$.

We claim $F<0$ at every point of $K$ where $S=0$. If all $r_i>0$,
this is the strict sign just proved. Otherwise, suppose exactly $\ell$
reciprocals vanish and delete the corresponding diverging eigenvalues
from an approximating sequence. Let $H_*$ be the limiting finite vector
and $P_*$ the product of its remaining $L_i$ entries. Deletion gives
$H_*\in\G_{k-\ell}$, and expansion by the diverging eigenvalues gives
\[
 S=\frac{\sigma_{k-\ell}(H_*)}{P_*},\qquad
 F=\frac{Q_{n-\ell,k-\ell,p-\ell}(H_*)}{P_*}.
\]
For example, the terms surviving in $\sigma_j/P$ are precisely those
containing all $\ell$ diverging entries; their remaining degree is
$j-\ell$. This also explains the shift $p\mapsto p-\ell$ in $Q$.
The selected eigenvalue of $H_*$ remains $1$ and
$k-\ell<p-\ell<k-\ell+1$. The strict sign in the smaller dimension
therefore proves the claim.

Since $S\geq0$ on $K$, compactness gives a number $\delta>0$ such
that $F>0$ implies $S\geq\delta$. Indeed, otherwise a sequence with
$F>0$ and $S\to0$ would have a limit in $K$ with $S=0$ and $F\geq0$,
contrary to the strict sign just proved. Also $F\leq C$ on $K$.
At a point where $Q_{n,k,p}(H)>0$,
\[
 Q_{n,k,p}(H)=PF\leq CP\leq CP^{p/k}
 \leq C\delta^{-p/k}(PS)^{p/k}
 =C\delta^{-p/k}\sigma_k(H)^{p/k}.
\]
If $Q\leq0$, the assertion is immediate. Homogeneity removes the
normalization $b=1$, and continuity gives the assertion at $b=0$.
\end{proof}

For a finite list $r$ and a noninteger real number $q$, set
\begin{equation}\label{eqV}
 V_q(r)=\sum_{j=0}^{\#r}\frac{(-1)^j\sigma_j(r)}{q-j}.
\end{equation}
For the empty list, $V_q(\varnothing)=1/q$.
In the next lemma lists are ordered decreasingly.
\begin{lemma}\label{lemfractional}
Let $r_1\geq\cdots\geq r_d$ belong to $\G_l$, where $1\leq l\leq d$,
and let $c=r_l>0$. If $l-1<q<l$, then
\[
 (-1)^{l+1}c^qV_q(r/c)\geq c^q\B(q-l+1,l)>0.
\]
The same lower bound holds with $c^qV_q(r/c)$ replaced by
$b^qV_q(r/b)$ whenever $b>0$ and
\begin{equation}\label{eqwindow}
 r_{l+1}\leq b\leq r_{l-1}.
\end{equation}
An endpoint whose index is outside $1,\ldots,d$ imposes no restriction.
\end{lemma}

\begin{proof}
Use $E_l$ and $T_l$ from Lemma~\ref{lemtaylor}.
Since $E_l(r,t)=O(t^l)$ at zero and $l-1<q<l$, termwise integration gives
\begin{align*}
 (-1)^lc^qV_q(r/c)
 &=-\int_{1/c}^\infty t^{-q-1}T_l(r,t)\,dt
   -\int_0^{1/c}t^{-q-1}(-1)^lE_l(r,t)\,dt\\
 &\leq-\int_{1/c}^\infty t^{-q-1}(ct-1)^{l-1}\,dt
 =-c^q\B(q-l+1,l).
\end{align*}
This proves the sign at $b=c$.

To see why it persists throughout \eqref{eqwindow}, define
$F(b)=(-1)^{l+1}b^qV_q(r/b)$ for $b>0$. Differentiation gives
\[
 F'(b)=(-1)^{l+1}b^{q-d-1}\prod_{i=1}^d(b-r_i).
\]
On the positive part of $[r_{l+1},r_l]$ this derivative is nonpositive,
because the product has $l$ negative factors. On
$[r_l,r_{l-1}]$ it is nonnegative, because there are $l-1$ such factors.
Thus $F$ attains its minimum over the stated interval at $c=r_l$.
Repeated roots give the same weak signs by continuity.
\end{proof}

\section{The Hessian obstacle argument}

We use the following obstacle and potential estimates.

First, if $h\in C^\infty(\overline{B_R})$ is positive on the boundary
and $\tau>0$, there is a $k$-admissible function
$z_\tau\in C^{1,1}(\overline{B_R})$ such that
\begin{equation}\label{eqobstacle}
 \begin{gathered}
 z_\tau=0\text{ on }\partial B_R,\qquad z_\tau\leq h,\\
 \sigma_k(D^2z_\tau)^{1/k}\geq\tau\text{ a.e. in }B_R,\\
 \sigma_k(D^2z_\tau)^{1/k}=\tau\text{ in }\{z_\tau<h\}.
 \end{gathered}
\end{equation}
This follows from \cite[Theorems~1.1 and~1.2]{JW14}, with
$\kappa=0$, boundary value zero, and right-hand side $\tau$.
The admissible subsolution is $C(|x|^2-R^2)$ for sufficiently large
$C$, and the auxiliary strictly convex function is $|x|^2$.
The operator $\sigma_k^{1/k}$ is among the examples covered by
their theorem. The $C^{1,1}$ viscosity inequalities hold almost
everywhere by testing with second-order Taylor polynomials plus
$\epsilon|x-x_0|^2$ and letting $\epsilon\downarrow0$.
On the noncontact set the solution is classical.

Second, for a zero-boundary $k$-admissible
$z\in C^{1,1}(\overline{B_R})$ and $p>\max\{k,n/2\}$,
\begin{equation}\label{eqheight}
 \norm z_{L^\infty(B_R)}
 \leq C_{n,k,p}R^{2-n/p}
       \norm{\sigma_k(D^2z)^{1/k}}_{L^p(B_R)}.
\end{equation}
Put $F=\sigma_k(D^2z)^{1/k}$. Subharmonicity gives $z\leq0$.
We first identify the Hessian measure used in the potential estimate.
On a relatively compact subdomain, mollify $z$ to obtain $z_\epsilon$.
The convexity of the elementary symmetric matrix cone implies
$D^2z_\epsilon\in\G_k$. Moreover, $z_\epsilon\to z$ locally uniformly,
so $z$ is also $k$-subharmonic in the viscosity sense; see
\cite[Section 2]{TrudingerWang}. The Hessians are locally uniformly
bounded and converge in every finite $L^s$ space. Hence
$\sigma_k(D^2z_\epsilon)\to F^k$ in $L^1$ locally. Weak continuity of
Hessian measures, in \cite[Theorem 7.1]{PV08}, gives
$\mu_k[z]=F^k\,dx$. Applying \cite[Theorem 7.4]{PV08} to
$-z$ on the ball now gives
\[
 -z(x)\leq C_{n,k}\int_0^{4R}
 \left(\frac{\int_{B_r(x)\cap B_R}F^k\,dy}{r^{n-2k}}\right)^{1/k}
 \frac{dr}{r}.
\]
The ball is uniformly $(k-1)$-convex, and
$F^k\in L^{p/k}(B_R)$ with $p/k>\max\{1,n/(2k)\}$,
so the hypotheses of the potential estimate hold.
H\"older's inequality bounds the expression in parentheses to the power
$1/k$ by
$C r^{2-n/p}\norm F_{L^p(B_R)}$.
The integral of $r^{2-n/p}\,dr/r$ over $(0,4R)$ is finite since $p>n/2$,
and equals a constant times $R^{2-n/p}$. This proves \eqref{eqheight}.

\begin{lemma}\label{lemobstacle}
Let $1\leq k\leq n-1$, $k<p<k+1$, and $p>n/2$.
Let $\mathcal A$ be a class of symmetric matrices uniformly bounded
above and below by positive multiples of the identity and set
$F_{\mathcal A}(H)=\inf_{A\in\mathcal A}A:H$.
Suppose a continuous function $Q$ on $\G_k$ has nonnegative integral
on every zero-boundary $k$-admissible $C^{1,1}$ function on a ball, and
there are positive constants $C_0,C_1,c_1$ such that
\begin{equation}\label{eqtwo-bounds}
 \begin{aligned}
 Q(H)&\leq C_0\sigma_k(H)^{p/k},\\
 Q(H)&\leq C_1\bigl(F_{\mathcal A}(H)^+\bigr)^p
                 -c_1\sigma_k(H)^{p/k}.
 \end{aligned}
\end{equation}
Then, for every measurable $A(x)\in\mathcal A$ and every
$u\in W^{2,p}(B_R)$,
\begin{equation}\label{eqball}
 \sup_{B_R}u\leq\sup_{\partial B_R}u^+
 +CR^{2-n/p}\norm{(A:D^2u)^-}_{L^p(B_R)}.
\end{equation}
\end{lemma}

\begin{proof}
First assume $u$ is smooth on the closed ball and $u\leq0$ on its
boundary. Let $\varepsilon>0$, set $h=-u+\varepsilon$, and choose
$z_\tau$ from \eqref{eqobstacle}. Denote its contact set by
$K=\{z_\tau=h\}$ and write $f=(A:D^2u)^-$.
At almost every point of $K$, the function $h-z_\tau$ has a local
minimum zero and is twice differentiable. Hence
\[
 D^2z_\tau\leq D^2h=-D^2u.
\]
Since $A$ is positive definite,
\[
 F_{\mathcal A}(D^2z_\tau)
 \leq A:D^2z_\tau\leq-A:D^2u\leq f.
\]
Use the second inequality of \eqref{eqtwo-bounds} on $K$.
On $B_R\setminus K$, use the first inequality and
$\sigma_k(D^2z_\tau)^{1/k}=\tau$. Nonnegativity of the integral of $Q$
then gives
\[
 0\leq\int_{B_R}Q(D^2z_\tau)\,dx
 \leq C_1\int_K f^p\,dx
       -c_1\int_K\sigma_k(D^2z_\tau)^{p/k}\,dx
       +C_0\tau^p|B_R\setminus K|.
\]
Adding the known contribution off $K$, we conclude that
\[
 \int_{B_R}\sigma_k(D^2z_\tau)^{p/k}\,dx
 \leq C\int_{B_R}f^p\,dx+C\tau^p|B_R|.
\]
Since $z_\tau\leq-u+\varepsilon$, we have $u-\varepsilon\leq-z_\tau$.
Applying \eqref{eqheight} yields
\[
 \sup_{B_R}u-\varepsilon
 \leq CR^{2-n/p}\norm f_{L^p(B_R)}+C\tau R^2.
\]
Let $\tau\downarrow0$ and then $\varepsilon\downarrow0$.
 Subtracting
$\sup_{\partial B_R}u^+$ gives \eqref{eqball} for arbitrary smooth $u$.

Finally approximate a function in $W^{2,p}(B_R)$ by functions smooth on
the closed ball. Since $p>n/2$, the approximation is also uniform on
$\overline{B_R}$. The coefficient field is fixed and bounded, so
$A:D^2u_j\to A:D^2u$ in $L^p$. The negative-part map is Lipschitz,
and therefore the right-hand sides converge as well. This proves
\eqref{eqball}.
\end{proof}

\section{The matrix inequality for fixed trace}
For the fixed-trace class, write $A=\delta I+B$, where $B\geq0$ and
$\tr B=1-n\delta$. Minimizing $B:H$ places its trace in an eigendirection
of the least eigenvalue of $H$, giving $F_\delta$ in \eqref{eqextremals}.
We will show that $Q_{n,k,p}(H)<0$ whenever $H\ne0$ and $F_\delta(H)\leq0$,
provided $p>n[1-(n-1)\delta]$.
Throughout this section, lists other than the original eigenvalues
$\lambda_i(H)$ are ordered decreasingly. We use $V_q$ from \eqref{eqV}.

The following polynomial allows us to treat a weighted sum of deleted
pairs at once; the eigenvalue separation used in its proof is the
standard min--max principle, as in \cite{HornJohnson}. Its individual summands need not have the same sign.

\begin{lemma}\label{lempair-polynomial}
Let $N\geq2$, let $\mu_1\geq\cdots\geq\mu_N$ be a nonconstant real
vector, and set $S=\sum_i\mu_i$. Define
\[  \begin{aligned}
 P(t) &=\prod_{i=1}^N(t-\mu_i),\\
 R(t) & =\sum_{i<j}(\mu_i-\mu_j)^2
                 \prod_{h\ne i,j}(t-\mu_h).
                  \end{aligned}
\]
Then
\[
 R(t)=(Nt-S)P'(t)-N^2P(t)
      =D_2\prod_{i=1}^{N-2}(t-r_i),
      \]
where  $D_2=N\sum_i\mu_i^2-S^2>0,$
 the roots are real and, in decreasing order, satisfy
\begin{equation}\label{eqinterlace}
 \mu_{i+2}\leq r_i\leq\mu_i\qquad(1\leq i\leq N-2).
\end{equation}
If $2\leq k\leq N-1$ and $\mu\in\Gamma_k$, then
$r\in\Gamma_{k-1}$.
\end{lemma}

\begin{proof}
For $t$ away from the entries of $\mu$, multiply out the rational
expression for $R/P$. One obtains
\[
 \frac{R(t)}{P(t)}
 =\sum_{i<j}\frac{(\mu_i-\mu_j)^2}{(t-\mu_i)(t-\mu_j)}
 =(Nt-S)\sum_i\frac1{t-\mu_i}-N^2.
\]
This proves the first identity by polynomial continuation.
The degree is $N-2$ and its leading coefficient is
$\sum_{i<j}(\mu_i-\mu_j)^2=D_2>0$.

To locate the roots, use a symmetric compression.
Let $D=\diag(\mu)$ and let $V$ be the $N$ by $2$
matrix with columns $(1,\ldots,1)^T$ and $\mu$. These columns are
independent. Let $W$ have orthonormal columns spanning their orthogonal
complement. The Cauchy--Binet formula gives
\[
 P(t)\det\bigl(V^T(tI-D)^{-1}V\bigr)=R(t).
\]
In fact, the $2$ by $2$ minor of $V$ in rows $i,j$ is
$\mu_j-\mu_i$, so each term is exactly the corresponding term of $R$.
Choose also an orthonormal basis $Z$ for the range of $V$ and write
$V=ZB$. The complementary-minor identity for the inverse of $tI-D$,
in the orthonormal basis $(Z,W)$, yields
\[
 P(t)\det\bigl(V^T(tI-D)^{-1}V\bigr)
 =\det(V^TV)\det(tI-W^TDW).
\]
These identities initially hold away from finitely many values of $t$
and therefore hold as polynomial identities. Since
$\det(V^TV)=D_2$, the roots $r_i$ are the eigenvalues of the real
symmetric matrix $W^TDW$. Restricting a symmetric quadratic form to a
subspace of codimension two gives
$\mu_{i+2}\leq r_i\leq\mu_i$ by the min--max characterization of
eigenvalues. For $N=2$, the compressed matrix has size zero and its
determinant is $1$.

It remains to prove the cone assertion. Put $e_j=\sigma_j(\mu)$.
The coefficient of $(-1)^{j-2}t^{N-j}$ in $R$ is
\begin{equation}\label{eqDj}
 D_j=(N-j+1)Se_{j-1}-Nj e_j,\qquad 2\leq j\leq N.
\end{equation}
We show that $D_j>0$ for $2\leq j\leq k+1$.
If $e_j\leq0$, then $j=k+1$, while $S>0$ and $e_{j-1}>0$,
so the assertion follows immediately. If $e_j>0$, then
$\mu\in\Gamma_j$. Newton's inequalities say that the positive ratios
\[
 \frac{e_h/\binom Nh}{e_{h-1}/\binom N{h-1}}
 \quad(1\leq h\leq j)
\]
are nonincreasing. The ratio for $h=2$ is strictly smaller than the
ratio for $h=1$, since
\[
 \left(\frac SN\right)^2-\frac{e_2}{\binom N2}
 =\frac{N\sum_i\mu_i^2-S^2}{N^2(N-1)}>0.
\]
Consequently
\[
 \frac{e_j}{e_{j-1}}<\frac{N-j+1}{Nj}S,
\]
which again gives $D_j>0$. Finally
$\sigma_h(r)=D_{h+2}/D_2>0$ for $1\leq h\leq k-1$.
\end{proof}

\begin{proposition}\label{proptrace-sign}
Let $1\leq k\leq n-1$ and $k<p<k+1$.
Suppose $H\in\G_k$ has a negative smallest eigenvalue and its $k$th
largest eigenvalue is positive. Divide all eigenvalues by the latter and
write the normalized vector as
\[
 H=\diag(-a,\mu_1,\ldots,\mu_N),\qquad
 N=n-1,\quad a>0,\quad \mu_1\geq\cdots\geq\mu_N,\quad\mu_k=1.\]
Set $S=\sum_i\mu_i$ and $a_0=(n-p)S/(pN)$.
If $a\geq a_0$, then $Q_{n,k,p}(H)\leq0$.
Equality holds precisely when $a=a_0$ and every $\mu_i=1$.
\end{proposition}

\begin{proof}
Write $e_j=\sigma_j(\mu)$. Since
$\sigma_j(H)=e_j-ae_{j-1}$, the recurrence
\[
 e_0=1,\qquad e_j=\sigma_j(H)+ae_{j-1}
\]
shows that $e_j>0$ for $1\leq j\leq k$.
Thus $\mu\in\Gamma_k$ and $S>0$. The definition of $Q$ gives
\begin{equation}\label{eqaffine}
 Q_{n,k,p}(H)=(-1)^k\bigl[V_p(\mu)+aV_{p-1}(\mu)\bigr].
\end{equation}
Apply Lemma~\ref{lemfractional} with $l=k$, $q=p-1$, and $c=1$.
It gives $(-1)^kV_{p-1}(\mu)<0$.
The right side of \eqref{eqaffine}, with $\mu$ fixed, is therefore
strictly decreasing as $a$ increases. It suffices to evaluate it at $a_0.$

We claim the exact identity
\begin{equation}\label{eqpair-identity}
 V_p(\mu)+a_0V_{p-1}(\mu)
 =-\frac1{pN}\sum_{i<j}(\mu_i-\mu_j)^2
                       V_{p-2}(\mu\setminus\{\mu_i,\mu_j\}).
\end{equation}
Here the $i$th and $j$th entries are omitted. To check the identity, expand the right-hand weighted
sum and compare coefficients in \eqref{eqDj}. It equals
\[
 \sum_{j=2}^N\frac{(-1)^j}{p-j}
       \bigl[(N-j+1)Se_{j-1}-Nj e_j\bigr].
\]
One may add the terms $j=0,1,N+1$ by setting $e_j=0$ outside
$0\leq j\leq N$; these added coefficients vanish.
Use $N-j+1=(p-j)+(n-p)$ and $j=p-(p-j)$.
The terms without a denominator cancel because
\[
 \sum_{j=0}^N(-1)^je_j=\prod_i(1-\mu_i)=0,
\]
and the remaining terms are
$-pN V_p(\mu)-(n-p)S V_{p-1}(\mu)$.
This proves \eqref{eqpair-identity}.

If $\mu$ is constant, all its entries equal $1$ and the weighted sum
vanishes. Suppose it is not constant and write $R(t)=D_2\prod(t-r_i)$
as in Lemma~\ref{lempair-polynomial}. Comparing the coefficients once
more shows that the sum in \eqref{eqpair-identity} is
\[
 D_2V_{p-2}(r).
\]
For $2\leq k\leq n-2$, the root vector belongs to $\Gamma_{k-1}$.
Moreover, \eqref{eqinterlace} and $\mu_k=1$ imply
\[
 r_k\leq1\leq r_{k-2},
\]
with nonexistent endpoints omitted. These are exactly the restrictions
\eqref{eqwindow} with $l=k-1$ and $b=1$.
The $(k-1)$st largest root is positive by the open-cone conclusion of
Lemma~\ref{lemcone}. Lemma~\ref{lemfractional}, with $q=p-2$, gives
$(-1)^kV_{p-2}(r)>0$. Combining this sign with
\eqref{eqaffine}--\eqref{eqpair-identity} proves strict negativity at
$a_0$.

The two end values of $k$ can be handled directly. If $k=1$, then
$p-2\in(-1,0)$ and every $r_i\leq\mu_1=1$.
Termwise integration gives
\[
 V_{p-2}(r)=-\int_0^1t^{1-p}\prod_i(1-tr_i)\,dt<0.
\]
The formula also holds for an empty root vector. If $k=n-1=N$, then
every $\mu_i\geq1$. Any vector $\nu$ obtained by deleting two entries
has length $n-3$ and satisfies
\[
 V_{p-2}(\nu)=\int_0^1t^{p-n}\prod_{i=1}^{n-3}(t-\nu_i)\,dt.
\]
This integral has strict sign $(-1)^{n-3}=(-1)^k$; it converges because
$p>n-1$. At least one pair has positive weight when $\mu$ is
nonconstant. The weighted sum therefore has the same strict sign.

We have proved strict negativity at $a_0$ in every nonconstant case,
and equality there in the constant case. Strict decrease in $a$ proves
the asserted inequality and equality characterization.
\end{proof}

The following elementary observation turns a strict matrix sign into
an estimate with uniform constants. We will use it for both coefficient
classes.

\begin{lemma}\label{lemcoercivity}
Let $\mathcal C$ be a closed cone in a finite-dimensional normed space.
Let $Q$ and $F$ be continuous on $\mathcal C$, positively homogeneous
of degrees $p>0$ and one, respectively. Suppose
\[
 H\in\mathcal C\setminus\{0\},\quad F(H)\leq0
 \quad\Longrightarrow\quad Q(H)<0.
\]
There are $C,c>0$ such that
\[
 Q(H)\leq C(F(H)^+)^p-c|H|^p\qquad(H\in\mathcal C).
\]
\end{lemma}
\begin{proof}
On the compact set $K=\mathcal C\cap\{|H|=1\}$, first suppose
$\{F\leq0\}$ is nonempty. Its strict negative maximum for $Q$,
and continuity, give $\delta,\varepsilon>0$ such that
$F\leq\varepsilon$ implies $Q\leq-\delta$ on $K$.
Choose $0<c<\delta$. On the remaining part of $K$, $F>\varepsilon$
and $Q$ is bounded above, so a sufficiently large $C$ gives
$Q+c\leq C(F^+)^p$. If $\{F\leq0\}$ is empty, $F$ has a positive
minimum on $K$ and the same estimate follows directly.
Homogeneity proves the assertion.
\end{proof}

\begin{proposition}\label{proptrace-coercivity}
Let $0<\delta\leq1/n$, $k<p<k+1$, and $p>n[1-(n-1)\delta]$.
There are constants $C,c>0$, depending only on $n,k,p,\delta$, such that
\begin{equation}\label{eqtrace-coercivity}
 Q_{n,k,p}(H)\leq C\bigl(F_\delta(H)^+\bigr)^p-c|H|^p
 \qquad(H\in\G_k).
\end{equation}
In particular the last term may be replaced by
$-c'\sigma_k(H)^{p/k}$ for a suitable $c'>0$.
\end{proposition}

\begin{proof}
First suppose $H\ne0$ and $F_\delta(H)\leq0$.
A nonzero positive semidefinite matrix has $F_\delta(H)>0$.
Thus $H$ has a negative smallest eigenvalue, and
Lemma~\ref{lemcone} implies that the selected eigenvalue $b$ is positive.
Normalize $b=1$ and use the notation of Proposition~\ref{proptrace-sign}.
The condition $F_\delta(H)\leq0$ becomes
\[
 0\geq F_\delta(H)=\delta S-[1-(n-1)\delta]a.
\]
Since $S>0$,
\[
 a\geq\frac{\delta S}{1-(n-1)\delta}
   >\frac{(n-p)S}{p(n-1)}=a_0.
\]
The strict inequality is equivalent to the assumed restriction on $p$.
Strict decrease in \eqref{eqaffine} therefore gives $Q_{n,k,p}(H)<0$.

Lemma~\ref{lemcoercivity}, applied on $\overline\Gamma_k$,
now gives \eqref{eqtrace-coercivity}.
Finally, $\sigma_k(H)^{p/k}\leq C_{n,k,p}|H|^p$.
\end{proof}

\section{The matrix inequality for general ellipticity bounds}
Fix $1\leq k<n$ and $k<p<k+1$. We prove the matrix sign for
$\Mm$ by reducing the spectrum to two distinct nonzero values.
The negative eigenvalues can be averaged. The positive eigenvalues need
a separate argument, since moving them also changes the selected
eigenvalue in the definition of $Q$.
The construction of $Q_{d,k,p}$ and the results of Section~4 apply in
each dimension $d\geq k+1$. Eigenvalues in this section are ordered
increasingly unless stated otherwise.

\begin{lemma}\label{lemnegative-average}
Suppose $H\in\G_k\subset\Sym(d)$ has a negative eigenvalue.
Increasing the magnitude of one negative eigenvalue strictly decreases
$Q_{d,k,p}(H)$ as long as the resulting matrix remains in $\G_k$.
Replacing two negative eigenvalues by their arithmetic mean preserves
$\G_k$, the selected eigenvalue, and the positive and negative spectral
sums. This replacement strictly increases $Q$ unless the entries are
equal. Consequently, replacing all negative entries by their common
mean does not decrease $Q$.
\end{lemma}
\begin{proof}
The selected eigenvalue is positive by Lemma~\ref{lemcone}; normalize
it to one. Write the spectrum as $(-a,\xi)$, where $a>0$.
Deleting the negative entry gives $\xi\in\Gamma_k$, and its $k$th
largest entry is still one. Multiaffinity of the elementary symmetric
functions gives
\[
 Q_{d,k,p}(-a,\xi)=(-1)^k\{V_p(\xi)+aV_{p-1}(\xi)\}.
\]
Lemma~\ref{lemfractional}, applied with index $k$ and exponent
$p-1$, shows that
\[
 \frac{\partial}{\partial a}Q_{d,k,p}(-a,\xi)
       =(-1)^kV_{p-1}(\xi)<0.
\]

For two negative entries write the spectrum as $(-a,-c,\xi)$ and
keep $s=a+c$ fixed. Deleting both entries gives $\xi\in\Gamma_k$,
with $k$th largest entry one. We have
\[
 \begin{aligned}
 \sigma_j(H)&=\sigma_j(\xi)-s\sigma_{j-1}(\xi)
                         +ac\sigma_{j-2}(\xi),\\
 Q_{d,k,p}(H)&=(-1)^k\{V_p(\xi)+sV_{p-1}(\xi)+acV_{p-2}(\xi)\}.
 \end{aligned}
\]
Symmetric functions with negative indices are zero. Increasing $ac$
preserves $\sigma_j(H)\geq0$ for $j\leq k$.
For $k\geq2$, Lemma~\ref{lemfractional} with index $k-1$,
exponent $p-2$, and cutoff one gives $(-1)^kV_{p-2}(\xi)>0$.
The cutoff is permitted because it is the $k$th largest entry of $\xi$.
For $k=1$, every entry of $\xi$ is at most one and $q=p-2\in(-1,0)$,
so
\[
 V_q(\xi)=-\int_0^1t^{-q-1}\prod_i(1-t\xi_i)\,dt<0.
\]
This gives the same positive coefficient of $ac$.
The arithmetic mean maximizes $ac$ at fixed $s$, proving the assertion
for two entries. Repeatedly averaging the least and greatest negative
magnitudes makes their variance tend to zero. Continuity proves the
assertion for all negative entries.
\end{proof}

\begin{lemma}\label{lempositive-extremum}
Fix a nonempty list of negative entries, a list of zero entries, and
the sum of the remaining positive entries. Allow the latter entries to
be nonnegative and require the full spectrum to belong to $\G_k$.
Suppose that this compact set is nonempty and that every global
maximizer of $Q_{d,k,p}$ belongs to $\Gamma_k$ and has every
variable entry positive.
Then some global maximizer has all variable entries equal.
\end{lemma}
\begin{proof}
The negative entries are fixed, so the selected eigenvalue $b$ stays
positive on the feasible set. Temporarily regard $b$ as an independent
variable and write
\[
 F(\lambda,b)=(-1)^k\sum_{j=0}^d
             \frac{(-1)^j\sigma_j(\lambda)b^{p-j}}{p-j}.
\]
This function is symmetric and affine in each eigenvalue separately.
Differentiation in $b$ gives
\[
 F_b(\lambda,b)=(-1)^kb^{p-d-1}\prod_{h=1}^d(b-\lambda_h).
\]
Thus $F_b(\lambda,b(\lambda))=0$ at the selected eigenvalue.
The selection $b(\lambda)$ is locally Lipschitz, and $F$ is smooth
where $b>0$. In a first-order expansion, the variation of the cutoff
therefore contributes only $O(|h|^2)$. It follows that
$Q(\lambda)=F(\lambda,b(\lambda))$ is $C^1$, including at repeated
selected eigenvalues, with $Q_i=F_i(\lambda,b(\lambda))$.

Assume that no maximizer has all variable entries equal. Choose a
maximizer where the multiplicity $r$ of $b$ is least.
The cone and positivity constraints are inactive at this point, so
$F_i=F_j$ for every pair of variable entries. Symmetry and multiaffinity
give the exact identity
\[
 F_i-F_j=(\lambda_j-\lambda_i)F_{ij}.
\]
For fixed $b$, the derivative $F_{ij}$ does not depend on
$\lambda_i$ or $\lambda_j$. Choose $i$ with $\lambda_i=b$ and
$j$ with $\lambda_j\ne b$. Then $F_{ij}=0$.

First suppose $r=1$. Under the sum-preserving variation
$\lambda_i(t)=b+t$ and $\lambda_j(t)=\lambda_j-t$, the selected
eigenvalue is $b(t)=b+t$ for small $|t|$. At $t=0$,
\[
 F_{bb}=(-1)^kb^{p-d-1}\prod_{h\ne i}(b-\lambda_h)<0,
 \qquad F_{ib}=-F_{bb},\qquad F_{jb}=0.
\]
There are exactly $k-1$ entries greater than the simple selected
eigenvalue, which proves the displayed sign. Since $F_{ii}=F_{jj}=0$,
\[
 \left.\frac{d^2}{dt^2}Q(\lambda(t))\right|_{t=0}
 =-2F_{ij}+2(F_{ib}-F_{jb})+F_{bb}=-F_{bb}>0.
\]
This contradicts maximality.

Now suppose $r\geq2$. Keep the cutoff fixed at $b$ in the same
variation. The change of $F$ is an exact quadratic polynomial, whose
linear and quadratic coefficients are $F_i-F_j$ and $-F_{ij}$.
Both vanish, so $F(\lambda(t),b)$ is constant.
To choose a sign of $t$ for which the selected eigenvalue remains $b$,
list the entries temporarily in decreasing order. Let the block equal
to $b$ occupy positions $L,\ldots,H$, where $L\leq k\leq H$.
If $k>L$, take $t>0$; if $k=L$, take $t<0$.
For sufficiently small nonzero $t$, the selected eigenvalue remains $b$
and has multiplicity $r-1$. The other changed entry remains separated
from $b$, and the cone and positivity conditions remain satisfied.
We have found a maximizer with smaller multiplicity, a contradiction.
\end{proof}

For later use put
\begin{equation}\label{eqtwo-value-spectrum}
 H_{d,m,a}=\diag(\underbrace{-a,\ldots,-a}_{m},
                      \underbrace{1,\ldots,1}_{d-m}),\qquad a>0.
\end{equation}
\begin{proposition}\label{proptwo-value-reduction}
Suppose that for every $k+1\leq d\leq n$, $1\leq m\leq d-k$, and
$a>0$ satisfying
\[
 H_{d,m,a}\in\G_k,\qquad (d-m)-\rho ma\leq0,
\]
one has $Q_{d,k,p}(H_{d,m,a})<0$.
Then $Q_{n,k,p}(H)<0$ for every nonzero $H\in\G_k$ with
$\Mm(H)\leq0$.
\end{proposition}
\begin{proof}
We use induction on the dimension, starting at $d=k+1$.
A nonzero positive semidefinite matrix has $\Mm(H)=\tr H>0$.
Thus every spectrum under consideration has a negative entry, a positive
selected eigenvalue, and at least $k$ positive entries.
The first part of the proof of Proposition~\ref{propupper} gives
$Q(H)<0$ when $\sigma_k(H)=0$. By Newton's inequalities, a point
of $\G_k$ with $\sigma_k>0$ is in $\Gamma_k$.

Deleting zero entries preserves cone membership, the selected
eigenvalue, both spectral sums, and $Q$. In fact, every symmetric
function of degree exceeding the new dimension vanishes, and all
remaining terms in \eqref{eqQ} are unchanged.
A face with a zero entry therefore reduces to a lower dimension with
the same $k,p,\rho$. The remaining dimension is at least $k+1$.
The induction hypothesis handles these faces; none occur in the base
dimension.

Suppose an interior-cone spectrum without zero entries satisfies
$Q\geq0$ and $\Mm\leq0$. Let $m$ be its number of negative entries.
Average them to $-a$ by Lemma~\ref{lemnegative-average}.
This does not decrease $Q$ and preserves both spectral sums.
Let $S$ be the positive sum. Maximize $Q$ over
\[
 \left\{x\in[0,S]^{d-m}:\sum_i x_i=S,
                         ((-a)^m,x)\in\G_k\right\}.
\]
The maximum is nonnegative, and all spectra in this set satisfy
$\Mm=S-\rho ma\leq0$. A maximizer cannot have $\sigma_k=0$
or a zero variable entry, by the preceding two paragraphs.
Lemma~\ref{lempositive-extremum} therefore gives a maximizer with
all positive entries equal to $t=S/(d-m)>0$. Homogeneity yields
\[
 Q_{d,k,p}((-a)^m,t^{d-m})
       =t^pQ_{d,k,p}(H_{d,m,a/t})<0,
\]
contrary to the nonnegative maximum.
\end{proof}

\section{Krawtchouk polynomials and the critical exponent}\label{secspectral}
The preceding reduction leaves a finite family of one-variable sign
questions. We identify their critical values explicitly.
For $d\geq2$, $1\leq m<d$, and $a>0$, define
\begin{equation}\label{eqRm}
 R_m(x,a)=m![z^m](1+az)^x(1-z)^{d-x}.
\end{equation}
Here $[z^m]$ denotes the Taylor coefficient at zero. The result is a
polynomial in $x$ of degree $m$, with leading coefficient $(1+a)^m$.
In the standard notation for Krawtchouk polynomials,
\[
 R_m(x,a)=(-1)^m m!\binom dm
      K_m\left(x;\frac1{1+a},d\right);
\]
see \cite[Eq.~18.23.3]{DLMF}. The generating-function identity first
holds for $x=0,\ldots,d$ and then for every $x$ by polynomial identity.

Let $K_{d,m}(a)$ be the tridiagonal matrix on $j=0,\ldots,m-1$ with
rates
\begin{equation}\label{eqrates}
 \lambda_j(a)=\frac{d-j}{1+a},\qquad
 \mu_j(a)=\frac{ja}{1+a}.
\end{equation}
Its diagonal, superdiagonal, and subdiagonal are respectively
$\lambda_j+\mu_j$, $-\lambda_j$, and $-\mu_j$.
At $a=(d-m)/(m\rho)$ it is $K_{d,m,\rho}$ from
\eqref{eqK-intro}.

\begin{lemma}\label{lemroots}
The roots of $R_m(\cdot,a)$ are simple and satisfy
$0<\kappa_1<\cdots<\kappa_m<d$. Between any two consecutive roots
there is an integer strictly between them. Moreover,
\begin{equation}\label{eqroot-coefficients}
 (-1)^mR_m(j,a)=\frac{m!\binom dm}{\binom dj}
                  \sigma_j(H_{d,m,a}),\qquad 0\leq j\leq d,
\end{equation}
and
\begin{equation}\label{eqroot-determinant}
 \det(K_{d,m}(a)-xI)=\frac{(-1)^mR_m(x,a)}{(1+a)^m}.
\end{equation}
The smallest eigenvalue of $K_{d,m}(a)$ is therefore positive and
simple and has a strictly positive right eigenvector.
\end{lemma}
\begin{proof}
The Krawtchouk polynomials are orthogonal on $\{0,\ldots,d\}$ with
positive weights
\[
 w_j=\binom dj\frac{a^{d-j}}{(1+a)^d}.
\]
The zeros of a degree-$m$ orthogonal polynomial are simple and lie
in the open convex hull of the support. See \cite[Sections~18.19 and~18.2(vi)]{DLMF}.
Separation by support points follows from orthogonality: if two consecutive
zeros had no support point strictly between them, division by the two
corresponding linear factors would give a polynomial of degree $m-2$
whose product with $R_m$ is nonnegative at every support point and
positive at some such point, contradicting orthogonality. The same
argument allows zeros at support points.

Changing $z$ to $-z$ in \eqref{eqRm} gives
\[
 \frac{(-1)^mR_m(j,a)}{m!}
 =\sum_q(-a)^q\binom jq\binom{d-j}{m-q}.
\]
Use the combinatorial identity
\[
 \binom dj\binom jq\binom{d-j}{m-q}
 =\binom dm\binom mq\binom{d-m}{j-q}.
\]
The sum on the right after this substitution is
$\binom dm\sigma_j(H_{d,m,a})$, proving
\eqref{eqroot-coefficients}.

For the determinant identity put $B=(1+a)x-d$ and
$F(z)=(1+az)^x(1-z)^{d-x}$.
Logarithmic differentiation gives
\[
 [1+(a-1)z-az^2]F'(z)=[B-adz]F(z).
\]
Comparison of coefficients yields $R_0=1$, $R_1=B$, and
\[
 R_j=[B-(a-1)(j-1)]R_{j-1}
       -a(j-1)(d-j+2)R_{j-2}\quad(j\geq2).
\]
The determinants of the leading blocks of
$(1+a)(K_{d,m}(a)-xI)$ satisfy this recurrence after multiplication
by $(-1)^j$, with initial values $1$ and $-B$.
This proves \eqref{eqroot-determinant}.

Finally let $\pi_j=\binom dj a^{-j}$. The identity
$\pi_j\lambda_j=\pi_{j+1}\mu_{j+1}$ shows that positive diagonal
similarity makes $K_{d,m}(a)$ symmetric with strictly negative
adjacent entries. Subtract this symmetric matrix from a sufficiently
large scalar multiple of the identity. The resulting nonnegative
irreducible matrix has a strictly positive eigenvector for its largest
eigenvalue. This is the eigenvector for the smallest eigenvalue of
$K_{d,m}(a)$, after reversing the similarity.
\end{proof}

\begin{proposition}\label{propspectral-sign}
For $1\leq k\leq d-1$,
\[
 H_{d,m,a}\in\G_k\quad\Longleftrightarrow\quad\kappa_1\geq k.
\]
If also $k<p<k+1$ and $k\leq d-m$, then under this cone condition
\begin{equation}\label{eqspectral-sign}
 \sgn Q_{d,k,p}(H_{d,m,a})=\sgn(\kappa_1-p).
\end{equation}
\end{proposition}
\begin{proof}
By \eqref{eqroot-coefficients}, the cone condition is equivalent to
$(-1)^mR_m(j,a)\geq0$ for $1\leq j\leq k$.
If $\kappa_1\geq k$, all factors $j-\kappa_i$ are nonpositive,
so these inequalities hold.
Conversely suppose $\kappa_1<k$. If $m=1$, or if $m\geq2$ and
$k<\kappa_2$, the sign at $j=k$ is wrong.
Otherwise $\kappa_2\leq k$, and Lemma~\ref{lemroots} gives an
integer strictly between $\kappa_1$ and $\kappa_2$. That integer is
at most $k$ and again gives the wrong sign.

Under the cone condition $\kappa_1\geq k$. If $m\geq2$, separation
by an integer implies $\kappa_2>k+1$. Indeed, if
$\kappa_1<k+1$, the separating integer is at least $k+1$;
if $\kappa_1\geq k+1$, the conclusion is immediate.
Thus only the first zero can lie in $(k,k+1)$.

Write $D_d(p)=\prod_{j=0}^d(p-j)$ and $\ell=d-m$.
The rational functions
\[
 \sum_{j=0}^d\frac{(-1)^j\sigma_j(H_{d,m,a})}{p-j},\qquad
 \frac{\ell!(-1)^\ell R_m(p,a)}{D_d(p)}
\]
have the same residues at every pole. This follows from
\eqref{eqroot-coefficients} and
$D_d'(j)=(-1)^{d-j}j!(d-j)!$.
Their difference has no poles and vanishes at infinity, so it is zero.
The selected eigenvalue is one, and hence
\[
 Q_{d,k,p}(H_{d,m,a})
 =\frac{(-1)^{k+d-m}(d-m)!R_m(p,a)}{D_d(p)}.
\]
On $(k,k+1)$, the sign of $D_d(p)$ is $(-1)^{d-k}$.
Consequently the sign of $Q$ is that of $(-1)^mR_m(p,a)$.
This is positive before the first zero, vanishes there, and is negative
between the first two zeros. For $m=1$ the last interval extends to
infinity. Formula \eqref{eqspectral-sign} follows.
\end{proof}

\begin{proposition}\label{propbox-coercivity}
Suppose $P_n(\rho)<p<n$, $p$ is not an integer, and $k=\lfloor p\rfloor$.
There are $C,c>0$, depending only on $n,p,\rho$, such that
\begin{equation}\label{eqbox-coercivity}
 Q_{n,k,p}(H)\leq C(\Mm(H)^+)^p-c|H|^p\qquad(H\in\G_k).
\end{equation}
\end{proposition}
\begin{proof}
Consider a spectrum $H_{d,m,a}$ from \eqref{eqtwo-value-spectrum} in
Proposition~\ref{proptwo-value-reduction}.
Its extremal pairing is nonpositive precisely when
$a\geq a_0=(d-m)/(m\rho)$.
Reducing the negative magnitudes to $a_0$ preserves cone membership.
The smallest root at $a_0$ is
$\kappa_{d,m}(\rho)\leq P_n(\rho)<p$, so
Proposition~\ref{propspectral-sign} gives
$Q_{d,k,p}(H_{d,m,a_0})<0$.
Increasing the magnitudes back to $a$ stays in the cone: each
intermediate spectrum is coordinatewise greater than or equal to
the final admissible spectrum. Lemma~\ref{lemnegative-average}
therefore gives $Q_{d,k,p}(H_{d,m,a})<0$.
Proposition~\ref{proptwo-value-reduction} now proves the strict sign
for every nonzero $H\in\G_k$ with $\Mm(H)\leq0$.
Apply Lemma~\ref{lemcoercivity} to obtain
\eqref{eqbox-coercivity}.
\end{proof}
\section{Proof of the upper estimates}\label{secupper}
We now have the three inequalities needed by
Lemma~\ref{lemobstacle}: the integral sign, the bound on the whole
cone, and the coercive bound for the chosen coefficient class.
For Theorem~\ref{Thm1.1}, choose a noninteger $q$ with
\[
 P_{\mathrm{tr}}(n,\delta)<q<\min\{p,n\},\qquad k=\lfloor q\rfloor.
\]
This is possible because the threshold  $P_{\mathrm{tr}}(n,\delta)$ is less than $n$.
Propositions~\ref{propintegral}, \ref{propupper}, and
\ref{proptrace-coercivity}, with $Q=Q_{n,k,q}$ and $F=F_\delta$,
give the hypotheses of Lemma~\ref{lemobstacle}.
Thus the ball estimate holds at exponent $q$.

For Theorem~\ref{Thm1.2}, Lemma~\ref{lemroots} gives
$0<\kappa_{d,m}(\rho)<d\leq n$, and the maximum is finite.
Hence $P_n(\rho)<n$. Choose a noninteger
\[
 P_n(\rho)<q<\min\{p,n\},\qquad k=\lfloor q\rfloor.
\]
Apply the same integral and whole-cone bounds, now using
Proposition~\ref{propbox-coercivity} and $F=\Mm$.
This proves the ball estimate for $I\leq A\leq\rho I$.
For $\lambda I\leq A\leq\Lambda I$, apply it to $A/\lambda$;
the forcing is divided by $\lambda$ as well.

In both cases H\"older's inequality gives
\[
 R^{2-n/q}\norm f_{L^q(B_R)}
 \leq |B_1|^{1/q-1/p}R^{2-n/p}\norm f_{L^p(B_R)}.
\]
For $u\in W^{2,p}(B_R)$, the inclusion in $W^{2,q}(B_R)$ therefore
proves the ball estimate at the requested exponent $p$.
This includes integers and $p\geq n$.
The intermediate $q$ is chosen using only the stated parameters,
so the constants have the required dependence.

It remains to pass from a ball to an arbitrary bounded open set while
assuming only local Sobolev regularity.
Let $D=\diam\Omega$, choose $x_0\in\Omega$, and set
\[
 M=\sup_{\partial\Omega}u^+,\qquad
 v=u-M-\varepsilon,\qquad \varepsilon>0.
\]
Then $\Omega\subset B_D(x_0)$.
Since $v\leq-\varepsilon$ on the boundary and $u$ is continuous on
$\overline\Omega$, the set $\{v\geq0\}$ is compactly contained in
$\Omega$.
Choose smooth convex functions $\phi_\delta:\R\to[0,\infty)$ with
\[
 \phi_\delta=0\text{ on }(-\infty,0],\quad
 0\leq\phi_\delta'\leq1,\quad \phi_\delta''\geq0,\quad
 \phi_\delta(t)\longrightarrow t^+
\]
uniformly on bounded intervals.
For example, integrate a smooth nondecreasing function which is zero
on $(-\infty,0]$ and one on $[\delta,\infty)$.
The chain rule gives
\[
 D^2\phi_\delta(v)=\phi_\delta'(v)D^2u
                    +\phi_\delta''(v)Du\otimes Du.
\]
Both terms are locally in $L^p$.
Indeed, $Du\in L^{2p}_{\mathrm{loc}}$: if $p<n$, Sobolev embedding
gives the larger exponent $np/(n-p)>2p$ because $p>n/2$;
if $p=n$, every finite target exponent is available; and if $p>n$,
the gradient is locally bounded.
Thus the zero extension of $w=\phi_\delta(v)$ belongs to
$W^{2,p}(B_D(x_0))$ and has compact support in $\Omega$.

Extend $A$ outside $\Omega$ by $I/n$ in the trace-one case and by
$\lambda I$ in the other case. Positivity of $A$ and convexity give
\[
 (A:D^2w)^-\leq\phi_\delta'(v)(A:D^2u)^-
                 \leq(A:D^2u)^-\quad\text{in }\Omega,
\]
and the left side vanishes outside $\Omega$.
The ball estimate for $w$, followed by $\delta\downarrow0$ and then
$\varepsilon\downarrow0$, proves \eqref{eqmaximum}.
This completes the upper estimates of Theorems~\ref{Thm1.1}
and~\ref{Thm1.2}.

\begin{corollary}\label{corvolume}
Assume the coefficient and exponent hypotheses of
Theorem~\ref{Thm1.1} or Theorem~\ref{Thm1.2}.
Let $\Omega$ be an open set of finite volume, and let $u$ belong
to the class \eqref{eqfunction-class}. If $u$ is bounded above, then
\begin{equation}\label{eqvolume-maximum}
 \sup_\Omega u\leq\sup_{\partial\Omega}u^+
 +C|\Omega|^{2/n-1/p}\norm{(A:D^2u)^-}_{L^p(\Omega)}.
\end{equation}
\end{corollary}

\begin{proof}
Suppose first that $p<\infty$. Put
$M=\sup_{\partial\Omega}u^+$ and $f=(A:D^2u)^-$.
For $\varepsilon>0$, use the convex functions $\phi_\delta$
from the preceding argument and extend
$w=\phi_\delta(u-M-\varepsilon)$ by zero.
Then $w\in W^{2,p}_{\mathrm{loc}}(\R^n)$, $w\geq0$,
and $(A:D^2w)^-\leq f$ on $\Omega$.
The zero extension is locally valid because $u-M-\varepsilon$
is negative near every boundary point. Extend $A$ outside $\Omega$
by a fixed matrix in the same coefficient class, and extend $f$ by zero.

Write $H=\sup w$. If $H>0$, choose $x_0$ with $w(x_0)>H/2$.
Take $0\leq\theta\leq1$ in $C_c^\infty(B_R(x_0))$, with
$\theta(x_0)=1$, $|D\theta|\leq C/R$, and
$|D^2\theta|\leq C/R^2$, and put $v=\theta^2w^2$.
Since $p>n/2$, the Sobolev chain rule gives
$v\in W^{2,p}(B_R(x_0))$ with zero boundary values.
Expanding and completing the square gives
\[
 \begin{aligned}
 A:D^2v
 &=2\theta^2w A:D^2w
   +2(\theta Dw+2wD\theta)^TA(\theta Dw+2wD\theta)\\
 &\quad+2\theta w^2 A:D^2\theta-6w^2(D\theta)^TAD\theta\\
 &\geq-2Hf-CR^{-2}H^2\ind_\Omega.
 \end{aligned}
\]
Applying the ball estimate yields
\[
 \frac{H^2}{4}
 \leq C R^{2-n/p}H\norm f_{L^p(\Omega)}
       +C R^{-n/p}|\Omega|^{1/p}H^2.
\]
Choose $R=K|\Omega|^{1/n}$ with $K$ sufficiently large to
absorb the last term. We obtain
$H\leq C|\Omega|^{2/n-1/p}\norm f_{L^p(\Omega)}$.
Letting $\delta\downarrow0$ and then $\varepsilon\downarrow0$
proves \eqref{eqvolume-maximum}.
For $p=\infty$, apply the finite-exponent estimate at a fixed
admissible exponent and use H\"older's inequality.
\end{proof}

\begin{corollary}\label{cordrift}
Assume the coefficient and exponent hypotheses of
Theorem~\ref{Thm1.1} or Theorem~\ref{Thm1.2}.
Let $\Omega$ have finite volume and let
$u\in W^{2,p}_{\mathrm{loc}}(\Omega)\cap C(\overline\Omega)$
be bounded above. Suppose
$\beta=b^TA^{-1}b\in L^p(\Omega)$ and
$f=(Lu)^-\in L^p(\Omega)$, where $L$ is given by \eqref{Operat-b}.
Let $C_0$ be the constant in \eqref{eqvolume-maximum} and set
\[
 E=\frac{C_0}{4}|\Omega|^{2/n-1/p}\norm\beta_{L^p(\Omega)}.
\]
Then
\[
 \sup_\Omega u\leq\sup_{\partial\Omega}u^+
 +C_0(e^{1+E}-1)|\Omega|^{2/n-1/p}\norm f_{L^p(\Omega)}.
\]
On a bounded open set, the same assertion holds with
$|\Omega|^{1/n}$ replaced by $\diam\Omega$ and $C_0$ by the
constant in \eqref{eqmaximum}.
\end{corollary}

\begin{proof}
We use the logarithmic change of variables from \cite{T26}.
Put $w=u-\sup_{\partial\Omega}u^+$ and $m=\sup_\Omega w$.
If $m\leq0$, the estimate is immediate. For $t>0$, set
\[
 v=\log\frac{m+t}{m+t-w}.
\]
The denominator is at least $t$. Since $p>n/2$,
$Du\in L^{2p}_{\mathrm{loc}}$, so the Sobolev chain rule gives
$v\in W^{2,p}_{\mathrm{loc}}(\Omega)\cap C(\overline\Omega)$.
Also $v\leq0$ on the boundary and
$\sup_\Omega v=\log(1+m/t)$.
Using $Lu\geq-f$ and completing the square, we obtain
\[
 A:D^2v\geq-\frac f{m+t-w}-b\cdot Dv+(Dv)^TADv
          \geq-\frac ft-\frac\beta4.
\]
The volume estimate applied to $v$ therefore gives
\[
 \log(1+m/t)
 \leq\frac{C_0|\Omega|^{2/n-1/p}}t\norm f_{L^p(\Omega)}+E.
\]
Take $t=m/(e^{1+E}-1)$ and rearrange.
The diameter version follows by applying \eqref{eqmaximum}
to the same function $v$.
\end{proof}

\section{Radial counterexamples}
The radial obstructions go back to Pucci \cite{Pucci}; see also
\cite{DK25}. We recall the cutoff constructions that give smooth
counterexamples at the critical exponents.

\begin{lemma}\label{lemlaplace-endpoint}
For $n\geq3$, there are nonnegative
$v_j\in C^\infty(\overline{B_1})$, vanishing on $\partial B_1$,
such that
\[
 \norm{v_j}_{L^\infty(B_1)}=v_j(0)=1,\qquad
 \norm{\Delta v_j}_{L^{n/2}(B_1)}\longrightarrow0.
\]
\end{lemma}
\begin{proof}
Choose a smooth radial cutoff $0\leq\chi_\varepsilon\leq1$, supported
in $(\varepsilon,e^{-2})$ and equal to one on
$(2\varepsilon,e^{-3})$, and put
\[
 f_\varepsilon(r)=\frac{\chi_\varepsilon(r)}{r^2\log(1/r)}.
\]
Its $L^{n/2}$ norm is uniformly bounded, because
\[
 \int_{B_1}f_\varepsilon^{n/2}\,dx
 \leq C_n\int_0^{e^{-2}}\frac{dr}{r[\log(1/r)]^{n/2}}
 =C_n\int_2^\infty t^{-n/2}\,dt<\infty.
\]
Let $w_\varepsilon$ be the radial solution of
$-\Delta w_\varepsilon=f_\varepsilon$ with zero boundary values.
It is constant near zero, smooth on the closed ball, nonnegative,
and decreasing in $r$. Integration of
$(r^{n-1}w_\varepsilon')'=-r^{n-1}f_\varepsilon$ gives
\[
 w_\varepsilon(0)=\frac1{n-2}\int_0^1
                (r-r^{n-1})f_\varepsilon(r)\,dr\longrightarrow\infty.
\]
The first term diverges as $\log\log(1/\varepsilon)$, whereas the
second is bounded. Set $v_\varepsilon=w_\varepsilon/w_\varepsilon(0)$
and choose any sequence $\varepsilon\downarrow0$.
\end{proof}

We now prove the necessity statement of Theorem~\ref{Thm1.1}.
If $\delta\geq1/[2(n-1)]$, then $P_{\mathrm{tr}}=n/2$, and
Lemma~\ref{lemlaplace-endpoint} applies with $A=I/n$.
Suppose $0<\delta<1/[2(n-1)]$, and set
\[
 c_\delta=1-(n-1)\delta,\qquad
 \mu=\frac{(n-1)\delta}{c_\delta}<1,\qquad
 p_c=nc_\delta=\frac n{1+\mu}>\frac n2.
\]
For $x\ne0$, let
\[
 A(x)=\delta I+(1-n\delta)e_r\otimes e_r,\qquad e_r=x/|x|.
\]
Its radial eigenvalue is $c_\delta$, its tangential eigenvalues are
$\delta$, and its trace is one. For nonnegative smooth radial $f$
vanishing near zero, the solution of $A:D^2w=-f$, with $w(1)=0$
and $w$ constant near zero, satisfies
\begin{equation}\label{eqradial-solution}
 \begin{aligned}
 w'(r)&=-c_\delta^{-1}r^{-\mu}\int_0^rs^\mu f(s)\,ds,\\
 w(0)&=\frac1{c_\delta(1-\mu)}
                     \int_0^1(r^\mu-r)f(r)\,dr.
 \end{aligned}
\end{equation}
The first identity is obtained by multiplying the radial equation by
$r^\mu$. Integrating $-w'$ and interchanging the two nonnegative
integrals gives the second.
Choose
\[
 f_\varepsilon(r)=\frac{\chi_\varepsilon(r)}
                         {r^{1+\mu}\log(1/r)}.
\]
Since $p_c(1+\mu)=n$ and $p_c>1$, its $L^{p_c}$ norm is uniformly
bounded. The $r^\mu f_\varepsilon$ term in
\eqref{eqradial-solution} diverges as $\log\log(1/\varepsilon)$;
the $rf_\varepsilon$ term remains bounded because $\mu<1$.
Normalize $w_\varepsilon$ by its maximum $w_\varepsilon(0)$.
The resulting functions have maximum one and forcing tending to zero
in $L^{p_c}$.

The coefficients can be made smooth without changing their action on
$w_\varepsilon$. This function is constant for $r\leq\varepsilon$.
Replace $A$ there by a smooth interpolation between $I/n$ for
$r\leq\varepsilon/2$ and the displayed radial matrix before
$r=\varepsilon$. Convexity preserves the trace and lower bound, and
the Hessian vanishes wherever $A$ is changed.
Finally, on the fixed ball, H\"older's inequality gives failure at
every $1\leq p\leq p_c$. This proves the lower assertion of
Theorem~\ref{Thm1.1}.

\section{Counterexamples at and below the spectral exponents}
A finite-order laminate on $\Sym(d)$ is a probability measure obtained
from a Dirac mass by finitely many rank-one splittings
\[
 \delta_H\longmapsto\theta\delta_{H_1}+(1-\theta)\delta_{H_2},
 \quad H=\theta H_1+(1-\theta)H_2,\quad
 \rank(H_1-H_2)=1,
\]
where $0<\theta<1$ and all matrices are symmetric. A splitting of an
atom of mass $\omega$ is multiplied by $\omega$.
We first state the Hessian realization needed to turn such a finite
measure into a smooth test function.

\begin{lemma}\label{lemrealization}
Let $d\geq2$, let $\nu=\sum_i\omega_i\delta_{H_i}$ be a finite-order
laminate on $\Sym(d)$ with barycenter $H$, and let
$G:\Sym(d)\to[0,\infty)$ be continuous.
For every $\varepsilon>0$ there is
$\phi\in C_c^\infty(B_1^d)$ such that
\[
 \norm\phi_{C^1}<\varepsilon,\qquad
 \left|\int_{B_1^d}G(H+D^2\phi)\,dx
               -|B_1^d|\int G\,d\nu\right|<\varepsilon.
\]
\end{lemma}
\begin{proof}
Apply the symmetric-matrix part of \cite[Proposition 2.3]{AFS} on an
inner ball $B_r^d$, $r<1$. It realizes the laminate by gradients of
Lipschitz maps, with affine boundary data, arbitrarily small uniform
distance from that affine map, and uniformly bounded gradient support
for the fixed laminate. The realizing map has symmetric gradient.
Since the ball is simply connected, it is the gradient of a function
$w\in W^{2,\infty}(B_r^d)$, with $Dw=Hx$ on the boundary.
The Hessian distribution approximates $\nu$.

Let $q=w-\tfrac12x\cdot Hx$. Then $Dq=0$ on the boundary.
Since the sphere is connected for $d\geq2$, its boundary value is
constant; subtract that constant. Now $q$ and $Dq$ both vanish on the
boundary. Smallness of $Dq$ also makes $q$ small in the uniform norm.
Extend $q$ by zero to $B_1^d$ and mollify at a scale smaller than the
distance between the two spheres. The extension belongs to
$W^{2,\infty}$, and the mollified functions are in $C_c^\infty(B_1^d)$.
They converge in $C^1$ and in every finite $W^{2,s}$ norm, with a
uniform Hessian bound for the fixed realization. Thus continuity of
$G$ passes its integral to the limit.

Choose $r$ close enough to one that the omitted shell has negligible
contribution, including its constant Hessian $H$. Then choose the
realization error small, and finally the mollification scale small.
These successive choices prove both stated inequalities.
\end{proof}

We also need smooth coefficients near the minimizing Pucci selection.
For a symmetric matrix $H$, define
\begin{equation}\label{eqsmooth-selection}
 A_\tau(H)=I+\frac{\rho-1}{2}
             \left[I-H(H^2+\tau^2I)^{-1/2}\right],\qquad \tau>0.
\end{equation}
This is smooth in $H$, satisfies $I\leq A_\tau(H)\leq\rho I$, and
\begin{equation}\label{eqselection-error}
 0\leq A_\tau(H):H-\Mm(H)\leq\frac{d(\rho-1)}2\tau.
\end{equation}
Diagonalization reduces these claims to
$-1<t/\sqrt{t^2+\tau^2}<1$ and
$0\leq |t|-t^2/\sqrt{t^2+\tau^2}\leq\tau$.
Thus small $L^p$ norm of $(\Mm(-D^2u))^+$ also gives small forcing
for smooth coefficients, by choosing $\tau$ small.

Fix $d\geq2$ and $1\leq m<d$. Write $\kappa(a)$ for the smallest
eigenvalue of $K_{d,m}(a)$ from \eqref{eqrates}.
For $h_m=0$, this matrix acts by
\[
 (K_{d,m}(a)h)_j=\frac{d-j}{1+a}(h_j-h_{j+1})
                 +\frac{ja}{1+a}(h_j-h_{j-1}),\quad 0\leq j<m,
\]
where the second term is absent at $j=0$.

\begin{lemma}\label{lemkappa-monotone}
The eigenvalue $\kappa(a)$ is positive and simple, and
$\kappa'(a)<0$. Its positive right eigenvector, normalized by $h_0=1$,
satisfies
\[
 1=h_0>h_1>\cdots>h_{m-1}>h_m=0.
\]
\end{lemma}
\begin{proof}
Positivity, simplicity, and the positive eigenvector follow from
Lemma~\ref{lemroots}. Let $\pi_j=\binom dj a^{-j}$.
Multiply the eigenvalue equations by $\pi_j$ and sum from $0$ to $j$.
Detailed balance $\pi_j\lambda_j=\pi_{j+1}\mu_{j+1}$ cancels the
interior terms and gives
\[
 \pi_j\lambda_j(h_j-h_{j+1})
          =\kappa\sum_{i=0}^j\pi_i h_i>0.
\]
This proves strict decrease of $h$.
A positive left eigenvector is $\ell_j=\pi_jh_j$.
The coefficients depend smoothly on $a$ and the eigenvalue is simple,
so differentiation of the eigenvalue equation yields
\[
 \kappa'(a)=\frac{\ell^TK_{d,m}'(a)h}{\ell^Th}.
\]
The derivative applied to $h$, with $h$ held fixed, has entries
\[
 (K_{d,m}'(a)h)_j
 =-\frac{d-j}{(1+a)^2}(h_j-h_{j+1})
      +\frac{j}{(1+a)^2}(h_j-h_{j-1})<0.
\]
The first term is strictly negative and the second is nonpositive.
The positive left eigenvector therefore gives $\kappa'(a)<0$.
\end{proof}

\begin{proposition}\label{proplaminate-counter}
Let $d\geq2$, $1\leq m<d$, and $\rho\geq1$.
If $1\leq p<\kappa_{d,m}(\rho)$, or if
$p=\kappa_{d,m}(\rho)>1$, then there are
$u_j\in C^\infty(\overline{B_1^d})$, zero on the boundary, and smooth
symmetric fields $I\leq A_j\leq\rho I$ such that
\[
 \norm{u_j}_{L^\infty(B_1^d)}\leq\frac34,\qquad
 \sup_{B_1^d}u_j\geq\frac14,\qquad
 \norm{(A_j:D^2u_j)^-}_{L^p(B_1^d)}\longrightarrow0.
\]
\end{proposition}
\begin{proof}
Set $a_0=(d-m)/(m\rho)$ and $a=a_0-\eta$, where
\begin{equation}\label{eqeta-range}
 0\leq\eta\leq\eta_0<a_0,
 \qquad \eta_0<\frac{d}{m\rho(m+1)}.
\end{equation}
At scale $t>0$, use diagonal matrices with entries in $\{t,-at\}$.
If a matrix has $j$ negative entries, then
\begin{equation}\label{eqstopped-force}
 \Mm(H)=t(d-j-\rho aj)
       =t\left[d\left(1-\frac jm\right)+\rho\eta j\right].
\end{equation}
We stop a branch at the end of a round when it has at least $m$
negative entries. If it has exactly $m$, its Pucci value is
$m\rho\eta t$. If it has more than $m$, its Pucci value is
nonpositive. To see the latter, the most restrictive condition is at
$j=m+1$, and it is exactly ensured by \eqref{eqeta-range}; the
quantity $(j-m)/j$ increases with $j$.

Choose $s\in(0,1)$. In one round, replace each of the $d$ coordinates,
one at a time, by either $t/s$ or $-at/s$.
The probability of the positive choice is
\[
 u_s=\frac{s+a}{1+a}\quad\hbox{for an old positive entry},\qquad
 z_s=\frac{a(1-s)}{1+a}\quad\hbox{for an old negative entry}.
\]
For example, $u_s(t/s)+(1-u_s)(-at/s)=t$; the other formula gives
$-at$. Thus each replacement preserves the parent matrix in expectation.
Only one diagonal entry changes, so every replacement is a symmetric
rank-one split. Stopping is decided only after all $d$ coordinates
have been replaced.

Let $T=T(s,a)$ be the transition matrix for the counts
$0,\ldots,m-1$ that have not been stopped. Independence of the coordinate
choices gives
\begin{equation}\label{eqtransition}
 T_{ij}=\sum_{r=\max\{0,j-d+i\}}^{\min\{i,j\}}
 \binom ir(1-z_s)^r z_s^{i-r}
 \binom{d-i}{j-r}(1-u_s)^{j-r}u_s^{d-i-j+r}.
\end{equation}
Here $r$ old negative entries remain negative, and $j-r$ old positive
entries become negative.
Write $\epsilon=1-s$. A positive entry changes sign with probability
$\epsilon/(1+a)$, while a negative entry changes sign with probability
$a\epsilon/(1+a)$. Events with at least two sign changes have probability
$O(\epsilon^2)$. Thus \eqref{eqtransition} gives, uniformly for $a\in[a_0-\eta_0,a_0]$,
\begin{equation}\label{eqtransition-expansion}
 T(s,a)=I-\epsilon K_{d,m}(a)+O(\epsilon^2).
\end{equation}
The loss of transitions to count $m$ is why the final diagonal of
$K$ includes the absorption rate.

Suppose $\kappa(a)>p$ and put $\Delta=\kappa(a)-p$.
The normalized positive eigenvector $h(a)$ is continuous in $a$.
After reducing $\eta_0$, there is $h_*>0$ such that
\[
 h_*\leq h_j(a)\leq1\quad(0\leq j<m).
\]
For each fixed $a$, choose $\epsilon>0$ small enough that the error
in \eqref{eqtransition-expansion}, applied to $h$, is at most
$(\Delta/2)\epsilon h$ componentwise. We may also require
\[
 Th\leq c h,\qquad c=1-(p+\Delta/2)\epsilon>0.
\]
Convexity of $s\mapsto s^p$ for $p\geq1$ gives
\[
 s^p-c\geq\frac{\Delta\epsilon}{2}>0.
\]
Fix such an $s$ for this $a$, and set $r=c/s^p<1$.

Starting at $I_d$, perform $N$ rounds on branches not yet stopped.
The resulting finite laminate $\nu_{\eta,N}$ has barycenter $I_d$.
Put $G(H)=(\Mm(H)^+)^p$.
Let $b_j$ be the probability that a round starting with $j<m$
negative entries ends with exactly $m$. This requires at least one
sign change, and each coordinate changes sign with probability at most
$\epsilon$. Therefore
\[
 0\leq b_j\leq d\epsilon,\qquad
 b\leq(d/h_*)\epsilon h.
\]
By \eqref{eqstopped-force}, a branch first stopped in round $q$ has scale $s^{-q}$ and positive Pucci value $m\rho\eta s^{-q}$ when its count is exactly $m$.
Writing $e_0$ for the first coordinate vector, its contribution summed
over all such branches is
\[
 \begin{aligned}
 E_{\mathrm{stopped}}
 &=(m\rho\eta)^p\sum_{q=1}^N s^{-qp}e_0^TT^{q-1}b\\
 &\leq C\eta^p\epsilon s^{-p}\sum_{q=0}^{N-1}r^q
 \leq\frac{C\eta^p\epsilon}{s^p-c}
 \leq\frac{C\eta^p}{\Delta}.
 \end{aligned}
\]
Branches stopped with more than $m$ negative entries have zero energy.
On a surviving branch the scale-one positive Pucci value is at most $d$.
Since $\mathbf1\leq h/h_*$ and $T^Nh\leq c^Nh$, the surviving
energy is at most $Cr^N$. Thus
\begin{equation}\label{eqlaminate-energy}
 \int G\,d\nu_{\eta,N}
       \leq\frac{C\eta^p}{\kappa(a)-p}+Cr^N.
\end{equation}
The constants are uniform in $\eta$, although $s$ and $r$ may depend
on $\eta$.

If $p<\kappa(a_0)$, take $\eta=0$, choose $s$ once, and let
$N\to\infty$. The first term in \eqref{eqlaminate-energy} vanishes.
If $p=\kappa(a_0)>1$, Lemma~\ref{lemkappa-monotone} gives
$\kappa(a_0-\eta)-p\geq c_0\eta$ for small positive $\eta$.
Consequently
\[
 \int G\,d\nu_{\eta,N}\leq C\eta^{p-1}+Cr^N.
\]
First let $\eta$ tend to zero. For each such $\eta$, choose $s$
as above and then a finite $N$ large enough to make $r^N$ tend to zero.
Here $\eta^{p-1}\to0$ because $p>1$.

Apply Lemma~\ref{lemrealization} separately to the resulting finite
laminates, with energy error tending to zero and
$\norm{\phi_j}_{C^0}<1/4$. Define
\[
 u_j(x)=\frac{1-|x|^2}{2}-\phi_j(x).
\]
These functions vanish on the boundary, have absolute value at most
$3/4$, and satisfy $u_j(0)>1/4$. Moreover,
$-D^2u_j=I_d+D^2\phi_j$, so $\int G(-D^2u_j)\to0$.
Take $A_j=A_{\tau_j}(-D^2u_j)$ from \eqref{eqsmooth-selection},
with $\tau_j\downarrow0$. By \eqref{eqselection-error},
\[
 (A_j:D^2u_j)^-
 =\bigl(A_j:(-D^2u_j)\bigr)^+
 \leq(\Mm(-D^2u_j))^++\frac{d(\rho-1)}2\tau_j.
\]
This proves the claimed forcing convergence with smooth coefficients.
\end{proof}

\section{Passage to higher dimensions}
We must account for each lower-dimensional term in
\eqref{eqbox-threshold}. A direct product cutoff would differentiate
the original test functions and give uncontrolled terms. The following
convex truncation avoids that problem and also preserves the endpoint.

\begin{lemma}\label{lemlifting}
Let $1\leq d<n$, $1\leq p<\infty$, and $\rho\geq1$.
Suppose $u_j\in C^\infty(\overline{B_1^d})$ vanish on the boundary
and there are measurable fields $I_d\leq A_j\leq\rho I_d$ and
constants $C,c>0$ such that
\[
 \norm{u_j}_\infty\leq C,\qquad \sup u_j\geq c,\qquad
 \norm{(A_j:D^2u_j)^-}_{L^p(B_1^d)}\longrightarrow0.
\]
Then there are $V_j\in C_c^\infty(B_1^n)$ and smooth symmetric
fields $I_n\leq B_j\leq\rho I_n$ such that
\[
 \norm{V_j}_\infty\leq C,\qquad \sup V_j\geq c/2,\qquad
 \norm{(B_j:D^2V_j)^-}_{L^p(B_1^n)}\longrightarrow0.
\]
\end{lemma}
\begin{proof}
Put $q=n-d$ and fix $0<\eta<c/4$.
Choose a smooth convex function $\chi:\R\to[0,\infty)$ with
\[
 \chi=0\text{ on }(-\infty,0],\quad
 0\leq\chi'\leq1,\quad\chi''\geq0,\quad
 ({t-\eta})^+\leq\chi(t)\leq t^+.
\]
For $0<\delta\leq1$, set
\[
 U_{j,\delta}(x,z)=u_j(x)-\delta|z|^2-\eta,\qquad
 v_{j,\delta}=\chi(U_{j,\delta}),\quad
 (x,z)\in B_1^d\times\R^q.
\]
For each $j$, continuity and the zero boundary values give a collar
of $\partial B_1^d$ where $u_j<\eta$. Hence extension by zero makes
$v_{j,\delta}$ smooth with compact support in $\R^n$, and
\[
 \supp v_{j,\delta}\subset B_1^d\times B_{\sqrt{C/\delta}}^q,
 \qquad \norm{v_{j,\delta}}_\infty\leq C,
 \qquad \sup v_{j,\delta}\geq c-2\eta>c/2.
\]
Extend $A_j$ by $I_d$ outside $B_1^d$ and put
$\widetilde A_j=\diag(A_j,I_q)$.
Writing $f_j=(A_j:D^2u_j)^-$, the chain rule and $\chi''\geq0$
give, on the support,
\[
 \begin{aligned}
 \widetilde A_j:D^2v_{j,\delta}
 &=\chi'(U_{j,\delta})(A_j:D^2u_j-2q\delta)
   +\chi''(U_{j,\delta})
       (DU_{j,\delta})^T\widetilde A_j DU_{j,\delta}\\
 &\geq-f_j-2q\delta.
 \end{aligned}
\]
The support bound therefore yields
\[
 \norm{(\widetilde A_j:D^2v_{j,\delta})^-}_{L^p(\R^n)}
 \leq C_1\left(\delta^{-q/(2p)}\norm{f_j}_{L^p(B_1^d)}
                       +\delta^{1-q/(2p)}\right).
\]
Set $R_\delta=2\sqrt{1+C/\delta}$ and define
$V_{j,\delta}(y)=v_{j,\delta}(R_\delta y)$ and
$B_{j,\delta}(y)=\widetilde A_j(R_\delta y)$.
The support lies in $B_{1/2}^n$. The Hessian gains a factor
$R_\delta^2$, while the $L^p$ norm gains $R_\delta^{-n/p}$.
Since $R_\delta$ is comparable to $\delta^{-1/2}$, we obtain
\begin{equation}\label{eqlift-error}
 \norm{(B_{j,\delta}:D^2V_{j,\delta})^-}_{L^p(B_1^n)}
 \leq C_2\left(\delta^{-1+d/(2p)}\norm{f_j}_{L^p(B_1^d)}
                              +\delta^{d/(2p)}\right).
\end{equation}
To make \eqref{eqlift-error} tend to zero, choose any sequence $\delta_\ell\downarrow0$. For each $\ell$, choose
$j_\ell$ large enough that the first term on the right tends to zero.
The second term tends to zero without a restriction on $p$.
This proves the assertions with measurable coefficients.

Finally let $H=-D^2V_{j_\ell,\delta_\ell}$. Since
$\Mm(H)\leq B_{j_\ell,\delta_\ell}:H$, the smooth selection
\eqref{eqsmooth-selection} has forcing no greater than the previous
forcing plus $n(\rho-1)\tau/2$.
Choose $\tau\downarrow0$ along the sequence to obtain the required
smooth fields.
\end{proof}

We can now complete Theorem~\ref{Thm1.2}.
For $1\leq p\leq n/2$, Lemma~\ref{lemlaplace-endpoint} gives
counterexamples with constant coefficients.
If $n/2<p<P_n(\rho)$, some pair $(d,m)$ in the finite maximum
\eqref{eqbox-threshold} satisfies $p<\kappa_{d,m}(\rho)$.
Proposition~\ref{proplaminate-counter} gives the counterexample in
dimension $d$; when $d<n$, apply Lemma~\ref{lemlifting}.
At $p=P_n(\rho)>n/2$, a maximizing pair satisfies
$p=\kappa_{d,m}(\rho)>1$, so the endpoint part of the same
proposition applies. The case $P_n(\rho)=n/2$ has already been treated.
Multiplying the coefficient fields by $\lambda$ restores
$\lambda I\leq A_j\leq\Lambda I$ and only multiplies their forcing
by $\lambda$. The necessity assertion is proved.

\end{document}